\documentclass[draft]{article}


\newcounter{parentnumber}

\usepackage{amsfonts}
\usepackage{amssymb}
\usepackage{textcomp}
\usepackage{latexsym,amsmath}
\usepackage{graphics}

\usepackage[utf8]{inputenc}
\usepackage{amsmath}
\usepackage{amsfonts}
\usepackage{amssymb}
\usepackage{calrsfs}

\usepackage{mathrsfs}
\usepackage{calrsfs}

\renewcommand{\phi}{\varphi}
\newcommand{\be}{\begin{equation}}
\newcommand{\ee}{\end{equation}}
\newcommand{\ba}{\begin{eqnarray}}
\newcommand{\ea}{\end{eqnarray}}
\newcommand{\ban}{\begin{eqnarray*}}
\newcommand{\ean}{\end{eqnarray*}}

\newcommand{\nul}{{\bf0}}
\newcommand{\rd}{{\mathbb R}^d}
\newcommand{\zd}{{\mathbb Z}^d}
\newcommand{\td}{{\mathbb T}^d}
\renewcommand{\r}{{\mathbb R}}
\newcommand{\z} {{\mathbb Z}}
\newcommand{\cn} {{\mathbb C}}
\newcommand{\n} {{\mathbb N}}

\newcommand{\ddd}{,\dots,}
\renewcommand{\lll}{\left(}
\newcommand{\rrr}{\right)}
\newcommand{\h}{\widehat}
\newcommand{\w}{\widetilde}

\def\N{{{\Bbb N}}}
\def\Z{{{\Bbb Z}}}
\def\T{{{\Bbb T}}}
\def\R{{\Bbb R}}

\def\vp{{\varphi}}
\def\t{{\theta }}

\def\C{{\Bbb C}}

\def\){\right)}
\def\({\left(}

\def\tt{\tilde{t}}
\def\kt{\tilde{k}}
\def\chit{\tilde{\chi}}

\def\supp{\operatorname{supp}}
\def\sinc{\operatorname{sinc}}
\def\mes{\operatorname{mes}}

\title{Quasi-projection  operators \\ in weighted $L_p$ spaces
\thanks{The first author (results of Sections 4.2 and 4.3 belong to this author) was partially supported by the DFG project KO 5804/1-1 and the project AFFMA that has received funding from the European Union's Horizon 2020 research and innovation programme under the Marie Sklodowska-Curie grant agreement No 704030; the second  author  (results of Sections 3 and 4.1 belong to this author) was supported by grant from Russian Science Foundation \# 18-11-00055.}}
\author{
Yu. Kolomoitsev$^{1}$ and M. Skopina$^{2}$
}
\date{\small $^{1}$University of L\"ubeck, Germany \\
\small 
 $^{2}$St. Petersburg State University, Russia \\
kolomoitsev@math.uni-luebeck.de, skopina@ms1167.spb.edu}

\begin{document}

\maketitle

\begin{abstract}
Approximation properties of  multivariate quasi-projection operators are studied.  Wide classes of such operators are considered, including
the sampling  and the Kantorovich-Kotelnikov type  operators generated by different band-limited functions.
The rate of convergence in the weighted $L_p$-spaces for these operators is investigated. The results allow us to  estimate the error for reconstruction of signals (approximated functions) whose decay is not enough to be in $L_p$.


\end{abstract}

\bigskip

\textbf{Keywords} Quasi-projection  operators, band-limited functions, approximation order, modulus of smoothness, matrix dilation, weighted $L_p$ spaces.

\medskip

\textbf{AMS Subject Classification}: 41A58, 41A25, 41A63


\newtheorem{theo}{Theorem}
\newtheorem{lem}[theo]{Lemma}
\newtheorem {prop} [theo] {Proposition}
\newtheorem {coro} [theo] {Corollary}
\newtheorem {defi} [theo] {Definition}
\newtheorem {rem} [theo] {Remark}
\newtheorem {ex} [theo] {Example}

\newtheorem{theorempart}{Theorem}[theo]

\newtheorem{lemmapart}{Lemma}[theo]

\newtheorem{proppart}{Proposition}[theo]

\newcommand{\tocsecindent}{\hspace{0mm}}

\section{Introduction}

The classical Kotelnikov formula (sampling expansion) provides exact reconstruction of band-limited signals based on the sampled values.
Up to now, an overwhelming diversity of digital signal processing applications and devices are based on it and more than
successfully use it. However, the class of band-limited signals is very small. 
To deal with sampling expansions for essentially wider classes of functions, one studies the convergence and error analysis of sampling expansions as the dilation factor goes to infinity (see, e.g., books~\cite{Stens}, \cite{Zayed}, and survey~\cite{Unser}).
In recent years, many works are dedicated to the study of approximation properties of  sampling expansions and their generalizations in $L_p$-norm (see~\cite{Butz4, Butz6, Butz5, JZ, KS1, KKS, Si1, Si2, Sk1}).

The classical sampling expansion $\sum_{k\in\z}f(M^{-j}k)\sinc(M^jx-k)$ is a special case of the
quasi-projection operators (or scaling expansions)
$$
Q_j(f; \w\vp, \vp)= \sum_{k\in\z} M^j \langle f,\w\vp(M^j\cdot-\,k)\rangle \vp(M^j\cdot-\,k)
$$
with   the Dirac
delta-function as $\w\vp$ and  the $\sinc$-function as $\vp$.

The operators  $Q_j(f; \w\vp, \vp)$ are actively
studied for different classes of functions $\vp$ and functions/distributions  $\w\vp$.
 In particular, these operators with some special functions $\vp$ and $\w\vp$ play an important role in wavelet theory, especially in the case of compactly supported $\vp$ and $\w\vp$. Approximation properties of  $Q_j(f; \w\vp, \vp)$ with compactly
supported $\vp$ and $\w\vp$ in the general situation were investigated  in~\cite{Jia1}, \cite{Jia2}, \cite{KS}.
Another important special case is the case where $\vp$ is a band-limited function and $\w\vp$ is locally  summable. Note that this class of  quasi-projection operators includes the classical Kantorovich--Kotelnikov operators, where $\w\vp$ is the characteristic
function of $[0,1]$. In this case $M^{j}\langle f,\w\vp(M^j\cdot-k)\rangle$ is
 the average value of $f$ near the node $M^{-j}k$.
  It is worth mentioning that using the averaging value instead of the exact  value $f(M^{-j}k)$ in the sampling expansion allows to deal with discontinuous signals and reduce the so-called time-jitter errors, which is very useful in Signal and Image Processing. Moreover, the Kantorovich--Kotelnikov type operators are bounded in $L_p(\R)$ and, therefore,  provide better approximation order than the sampling operators.
During recent years, the Kantorovich--Kotelnikov type  operators as well as their different generalizations and refinements have been especially actively studied (see, e.g.,~\cite{BDR, v58, Jia2,  KKS, KS, KS1, LJC, OT15, VZ2}).

In the multidimensional case, the quasi-projection operator $Q_j(f; \w\vp, \vp)$ can be defined using a matrix $M$ as dilation as follows
$$
Q_j(f; \w\vp, \vp)= \sum_{k\in\zd} |\det M|^{j} \langle f,\w\vp(M^j\cdot-\,k)\rangle \vp(M^j\cdot-\,k).
$$
For wide classes of  functions $\w\vp$ and
band-limited functions $\vp$, the multivariate quasi-projection operators with a matrix dilation
were studied in~\cite{KoS}.  Under the assumption that the Fourier transform
of  both the functions is $n+d+1$ times
continuously differentiable in a neighborhood of the origin
(but not necessarily continuous on $\rd$), and all derivatives up to
order $n-1$ of the function $1-\overline{\h\vp}\h{\w\vp}$ vanish at the origin, the following estimate was derived:
\begin{equation}
\label{100-}
  \|f-Q_j(f; \w\vp, \vp)\|_{p}\le C \omega_n\(f, \|M^{-j}\|\)_{p},\quad f\in L_p(\rd),\quad 1\le p\le\infty,
\end{equation}
where $\omega_n$ is the modulus of smoothness of order $n$ and $C$ does not depend on $f$ and $j$.

Now let us consider the case of sampling expansions.
If $Q_j(f; \w\vp, \vp)$ is the classical sampling expansion (i.e. $\vp$ is the $\sinc$-function and $\w\vp$ is the Dirac delta-function), then the well-known Brown's inequality~\cite{Brown}
$$
\left\|f- Q_j(f; \w\vp, \vp)\right\|_\infty\le
C\int\limits_{|\xi|>M^{j}/2}
 |\h f(\xi)|d\xi
$$
holds for every $j\in \Z_+$ and a function $f:\R\to \C$ whose Fourier transform  is summable on $\r$.
This Brown's result was strengthened in~\cite{KKS} in several directions. Namely,
a similar inequality was established in $L_p$-norm, $2\le p\le \infty$,
and the multivariate  operators $Q_j(f; \w\vp, \vp)$  with a matrix dilation were considered for a wide class of tempered distributions  $\w\vp$ and a wide class of band-limited functions $\vp$  instead of the $\sinc$-function.
In particular, in the case of the  Dirac delta-function as $\w\vp$, under the assumption that  the Fourier transform of $\vp$
is identical 1 on a $\delta$-neighborhood of the origin, for any
$\gamma>d/p$ and every function   $f$ such that  $\h f\in L_q(\R^d)$,   $\h f(\xi)=\mathcal{O}(|\xi|^{-\gamma-d/q})$,
$1/p+1/q=1$, the following inequality was derived:
	\be
	\|f-Q_j(f; \w\vp, \vp)\|_p\le C
	 \|M^{*-j}\| \(\,\,\int\limits_{|M^{*-j}\xi|\ge\delta}
	|\xi|^{q\gamma}| \h f(\xi)|^q d\xi\)^{1/q},
	 \label{8-}
	\ee
where  $C$ does not depend on $f$ and $j$.

Error estimates~\eqref{100-} and \eqref{8-} are aimed at the recovery of signals $f$,
but they are not applicable to  non-decaying signals and even for signals whose decay is not enough to belong to the space $L_p(\R^d)$.
The goal of the present paper is to obtain counterparts of these estimates in norms of some weighted $L_p$-spaces.
An idea to extend sampling theory in this regard was recently suggested by H.~Q.~Nguyen and M.~Unser~\cite{NU}.
In particular, they extended a number of basic facts of the theory of shift-invariant spaces in $L_p(\R^d)$ to
the weighted spaces $L_{p,1/w}(\R^d)$, where the function $w$ belongs to the class of the so-called submultiplicative weights. Since this class contains weights with polynomial growth,  signals which grow not faster than a polynomial are in
an appropriate space $L_{p,1/w}(\R^d)$.

Using the technique developed in~\cite{KoS} and \cite{KS1} as well as  several basic facts obtained in~\cite{NU} for the same class of submultiplicative weights,  we derive an analog of estimate~\eqref{100-} given in terms of the best approximation and moduli of smoothness in the weighted $L_p$ spaces (see Theorems~\ref{lemJack} and~\ref{th2}). However,  this technique does not allow to work with
slowly decaying functions $\phi$, in particular, with the $\sinc$-function. Thus, in contrast to the non-weighted case,
the corresponding  results  for $L_{p,1/w}$  are proved under  additional assumptions on the smoothness  of the Fourier transform of $\vp$.
To fix this drawback, we apply and extend  the technique developed in~\cite{Sk1} and \cite{KKS}. This is done in Sections~\ref{S5.2} and~\ref{S5.3}
but only for weights $w$ satisfying some additional properties.
Namely, the boundedness of the maximal function of the Hilbert transform in $L_{p,1/w}$ is required. For this purpose, we consider a certain subclass of Muckenhoupt type weights $\mathcal{A}_p(\R^d)$, which consists of all admissible weights $w$ that belong to the classical Muckenhoupt class $A_p(\R^1)$ in each variable uniformly with respect to other variables. Note that because of this we have to restrict  our consideration to diagonal dilation matrices $M$.

Approximation properties of the sampling expansions in $L_{p,1/w}$-norm for the submultiplicative
Muckenhoupt  weights $w$ and diagonal matrices $M$ are studied in Section~4.3.
In particular, under the additional assumption that $w$ is  band-limited and under certain assumptions on $\vp$, we obtain
an analog of estimate~(\ref{8-})  (see Theorem~\ref{theoQj_new}). Without this additional assumption on $w$ and for a larger class
of band-limited functions $\vp$, an error estimate for the sampling expansions is derived in Theorem~\ref{theoQj_new1}. Namely,
the Fourier transform of $\vp$ should be   smooth enough near the origin and  all derivatives up to
order $n-1$ of the function $1-{\h\vp}$ should vanish at the origin (instead of identical 1 in a neighborhood of the origin, as we suppose in Theorem~\ref{theoQj_new}). In this case, the approximation order represents a combination of estimates given in~\eqref{100-} and~\eqref{8-}.
As a consequence, we obtain the following approximation order:
$$
\bigg\Vert f- \sum\limits_{k\in\zd} f(M^{-j}k) \vp(M^j\cdot-\,k)\bigg\Vert_{p, 1/w}=
 \begin{cases}
	 	\mathcal{O}(|\lambda|^{-j(d/p+a)})  &\mbox{if } n>d/p+a,
\\
	 	\mathcal{O}(|\lambda|^{-jn}j^{1/2})  &\mbox{if } n=d/p+a,
\\
\mathcal{O}(|\lambda|^{-jn})  &\mbox{if }  n<d/p+a
	 	 \\
	\end{cases}
$$
whenever the Fourier
transform of the function $f/w$ belongs to $L_q(\R^d)$, $1/p+1/q=1$, and has decay of order $d+a$, $a>0$  (here $\lambda$ is the smallest (in absolute value) diagonal element of $M$).  Note that this estimate as well as Theorem~\ref{theoQj_new1} are new even in the non-weighted case.

The paper is organized as follows: in Section~2, we introduce notation and give  some basic facts. In Section~3,
we state and prove several auxiliary results.  Section~4 is dedicated to the main results.
In Subsection~4.1, we study approximation properties of the quasi-projection operator
$Q_j(f; \w\vp, \vp)$ in the case of $\vp$ and $\w\vp$ belonging to
certain classes of functions whose Fourier transform has some smoothness.
Only integrable functions  belong to this class. In Subsection 4.2, we solve similar problems
for the case of $\vp$ belonging to a class of band-limited functions  that includes non-integrable functions.
Subsection~4.3 is dedicated to studying sampling expansions (i.e., the case of the Dirac delta-function as $\w\vp$)
for a wide class of band-limited functions $\vp$.

\section{Notation and basic facts}
\label{notation}

$\n$ is the set of positive integers,
    $\r$  is the set of real numbers,
    $\cn$ is the set of complex numbers.
    $\rd$ is the
    $d$-dimensional Euclidean space,  $x = (x_1\ddd x_d)$ and $y =
    (y_1\ddd y_d)$ are its elements (vectors),
      $(x, y)~=~x_1y_1+~\dots~+x_dy_d$,
    $|x| = \sqrt {(x, x)}$, ${\bf0}=(0\ddd 0)\in\rd$;
      $B_r=\{x\in\rd:\ |x|\le r\}$, $\td=[-\frac 12,\frac 12]^d$;
    $\zd$ is the integer lattice
    in $\rd$, $\z_+^d:=\{x\in\zd:~x\geq~{\bf0}\}.$
    If $\alpha,\beta\in\zd_+$, $a,b\in\rd$, we set
    $[\alpha]=\sum\limits_{j=1}^d \alpha_j$,
    $\alpha!=\prod\limits_{j=1}^d(\alpha_j!),$
    $$\binom{\beta}{\alpha}=\frac{\beta!}{\alpha!(\beta-\alpha)!},\quad
       D^{\alpha}f=\frac{\partial^{[\alpha]} f}{\partial x^{\alpha}}=\frac{\partial^{[\alpha]} f}{\partial^{\alpha_1}x_1\dots
    \partial^{\alpha_d}x_d}.$$

A $d\times d$  real matrix $M$ whose
eigenvalues are bigger than 1 in modulus is called a  dilation matrix.
Since the spectrum of the operator $M^{-1}$ is
located in  $B_r$,
where $r=r(M^{-1}):=\lim_{j\to+\infty}\|M^{-j}\|^{1/j}$
is the spectral radius of $M^{-1}$, and there exists at least
one point of the spectrum on the boundary of $B_r$, we have
\be
	\|M^{-j}\|\le {C_{M,\vartheta}}\,\vartheta^{-j},\quad j\in \Z_+,
	\label{00++}
\ee
for every  positive number $\vartheta$  which is smaller in modulus
than any eigenvalue of  $M$.
In particular, we can take $\vartheta > 1$, then
$$
	\lim_{j\to+\infty}\|M^{-j}\|=0.
$$

In what follows,  the class of matrix dilations is denoted by $\frak M$, $m=|{\rm det}\,M|$,
$M^*$ is the conjugate {\sl transpose} of $M$, and  the $d\times d$  identity matrix is denoted by $I_d$,
note also that $M^0:=I_d$.

We will say that a weight $w\,:\ \R^d\mapsto [1,\infty)$  belongs to the class $\mathcal{W}^\alpha$ for some $\alpha>0$ if the following conditions hold:

\begin{description}
  \item[$1)$] $w$ is continuous and even;

  \item[$2)$]  there exists an even function  $w^* :\R^d\mapsto \R_+$ and a constant $c_w>0$ such that
\begin{equation}\label{W}
  w(x+y)\le  w^*(x)w(y),\quad w^*(x)\le c_w (1+|x|^2)^{\alpha/2}
\end{equation}
for all $x,y\in \R^d$;
  \item[$3)$] for every $M\in \mathfrak{M}$, there exists a constant $C'>0$ such that
  $$
  w^*(M^{-j}x)\le C' w^*(x)
  $$
for all  $x\in \rd$ and $j\in\z_+$.
		
 A model example of such a weight is
$$
w_\alpha(x)=(1+|x|^2)^{\alpha/2}, \ \ \alpha>0.
$$

\end{description}

Note that if $w\in \mathcal{W}^\alpha$, then condition~2)
is satisfied also  for $1/w$. Indeed,
\be
\label{1010}
\frac1{w(x+y)}=\frac{w(x)}{w(x+y)w(x)}
\le \frac{w^*(y)}{w(x)}.
\ee

Below, $L_{p,w}$ denotes the weighted space $L_{p}(\rd,w)$, $1\le p\le\infty$, with the norm $\Vert f\Vert_{p,w}=\Vert fw\Vert_{L_p(\R^d)}$.  In the unweighted case, we denote $L_p=L_p(\R^d)$. Obviously, if $w\in \mathcal{W}^\alpha$, then
$L_{p,w}\subset{L_p}$.
In  what follows, we will often use the following simple inequality for $f\in L_{p,1/w}$ and $w\in \mathcal{W}^\alpha$:
\begin{equation}\label{fwa}
  \Vert f(\cdot + h)\Vert_{p,1/w}\le  w^*(h)\Vert f\Vert_{p,1/w}.
\end{equation}
In particular, if $|h|\le 1$, then
\begin{equation}\label{fwah}
  \Vert f(\cdot + h)\Vert_{p,1/w}\le c_w 2^{\alpha/2} \Vert f\Vert_{p,1/w}.
\end{equation}

Similarly,  $\ell_{p,w}$ denotes the set of sequences $c=\{c_k\}_{k\in \Z^d}$
with the norm
$$
\Vert c\Vert_{\ell_{p,w}}=\(\sum_{k\in\Z^d} |c_kw(k)|^p\)^{1/p}.
$$

We use $W_{p,w}^n$, $1\le p\le\infty$, $n\in\n$,
to denote  the Sobolev space on~$\rd$, i.e. the set of
functions whose partial derivatives up to order $n$ are in $L_{p,w}$,
with the usual Sobolev semi-norm given by
$$
\|f\|_{\dot W_{p,w}^n}=\sum_{[\nu]=n}\Vert D^\nu f\Vert_{p,w}.
$$

If $f, g$ are functions defined on $\rd$ and $f\overline g\in L_1$,
then  $\langle  f, g\rangle:=\int_{\rd}f\overline g$.
If $f\in L_1$,  then its Fourier transform is $\mathcal{F}f(\xi)=\widehat
f(\xi)=\int_{\rd} f(x)e^{-2\pi i
(x,\xi)}\,dx$.

If $\phi$ is a function defined on $\rd$ and $M\in \mathfrak{M}$, we set
$$
\phi_{jk}(x):=m^{j/2}\phi(M^jx+k),\quad j\in\z,\,\, k\in\rd.
$$

Denote by $\mathcal{S}$ the Schwartz class of functions defined on $\rd$.
    The dual space of $\mathcal{S}$ is $\mathcal{S}'$, i.e. $\mathcal{S}'$ is
    the space of tempered distributions.
    The basic facts from distribution theory
    can be found, e.g., in~\cite{Vladimirov-1}.
    Suppose $f\in \mathcal{S}$, $\phi \in \mathcal{S}'$, then
    $\langle \phi, f\rangle:= \overline{\langle f, \phi\rangle}:=\phi(f)$.
    If  $\phi\in \mathcal{S}',$  then $\h \phi=\mathcal{F}\phi$ denotes its  Fourier transform
    defined by $\langle \h f, \h \phi\rangle=\langle f, \phi\rangle$,
    $f\in \mathcal{S}$.

    If  $\phi\in \mathcal{S}'$, $j\in\z, k\in\zd$, then we define $\phi_{jk}$ by
        $
        \langle f, \phi_{jk}\rangle=
        \langle f_{-j,-M^{-j}k},\phi\rangle$ for all $f\in \mathcal{S}$.

Denote by $\mathcal{S}_N'$ the set of all tempered distributions
	whose Fourier transform $\h{\phi}$ is a function on $\rd$
	such that $|\h{\phi}(\xi)|\le C_{\phi} |\xi|^{N}$
	 for almost all $\xi\notin\td$, $N=N({\phi})\ge 0,$ and
	 $|\h{\phi}(\xi)|\le C'_{{\phi}}$
	 for almost all $\xi\in\td$.

Denote by ${\cal B}={\cal B}(\R^d)$  the class of functions $\phi$ given by
\be
\phi(x)=\int\limits_{\rd}\theta(\xi)e^{2\pi i(x,\xi)}\,d\xi,
\label{zzz61}
\ee
where $\theta$ is supported in a parallelepiped $\Pi:=[a_1, b_1]\times\dots\times[a_d, b_d]$ and such that	$\theta\big|_\Pi\in C^d(\Pi)$.

Let $1\le p \le \infty$. Denote by ${\cal L}_{p}$ the set
	$$
	{\cal L}_{p}:=
	\left\{
	\phi\in L_{p}\,:\, \|\phi\|_{{\cal L}_{p}}:=
	\bigg\|\sum_{k\in\zd} \left|\phi(\cdot+k)\right|\bigg\|_{L_{p}(\td)}<\infty
	\right\}.
	$$
	With the norm $\|\cdot\|_{{\cal L}_p}$, ${\cal L}_p$ is a Banach space.
	The simple properties are:
	${\cal L}_1=L_1,$
	$\|\phi\|_p\le \|\phi\|_{{\cal L}_p}$,
	$\|\phi\|_{{\cal L}_q}\le \|\phi\|_{{\cal L}_p}$
	for $1\le q \le p \le\infty.$ Therefore, ${\cal L}_p\subset L_p$
	and ${\cal L}_p\subset {\cal L}_q$ for $1\le q \le p \le\infty.$
	If $\phi\in L_p$ and compactly supported, then $\phi\in {\cal L}_p$ for any $p\ge1.$
	If $\phi$ decays fast enough, i.e. there exist constants $C>0$
	and $\varepsilon>0$ such that
	$|\phi(x)|\le C( 1+|x|)^{-d-\varepsilon}$ for all $x\in\rd,$
	then $\phi\in {\cal L}_\infty$.

The corresponding weighted $\mathcal{L}_{p,w}$ space with respect to a weight function $w$ is defined according to the following weighted norm
$$
\Vert \vp\Vert_{{\cal L}_{p,w}}:=\Vert \vp w\Vert_{{\cal L}_{p}}.
$$

The modulus of smoothness $\omega_n(f,\cdot)_{{p,w}}$  of order $n\in \N$ for a function $f\in L_{p,w}$  is defined by
\begin{equation}\label{ClassMod}
  \omega_n(f,h)_{{p,w}}=\sup_{|\delta|<h,\, \delta\in \R^d} \Vert \Delta_\delta^n f\Vert_{{p,w}},
\end{equation}
where
$$
\Delta_\delta^n f(x)=\sum_{\nu=0}^n (-1)^\nu \binom{n}{\nu} f\(x+\delta \nu\).
$$

The modulus $\omega_n(f,\cdot)_{{p,w}}$ is the classical modulus of smoothness.
Together with the modulus~\eqref{ClassMod}, we will also use the following so-called anisotropic modulus of smoothness, in which the step $h$ is replaced by a $d \times d$ real matrix $M$:
$$
\Omega_n(f,M^{-1})_{{p,w}}=\sup_{|M\delta|<1,\, \delta\in \R^d} \Vert \Delta_\delta^n f\Vert_{{p,w}}.
$$

As usual, in the unweighted case, i.e., $w_0(x)\equiv 1$, we use the following notation: $\omega_n(f,h)_{{p}}=\omega_n(f,h)_{{p,w_0}}$ and $\Omega_n(f,M^{-1})_{{p}}=\Omega_n(f,M^{-1})_{{p,w_0}}$.

It is obvious that for any $M\in \mathfrak{M}$ and $f\in L_{p,w}$, one has
$$
\Omega_n(f,M^{-1})_{{p,w}}\le \omega_n(f,\Vert M^{-1}\Vert)_{{p,w}}.
$$
Note also that  in the case of the diagonal matrix $M$ with $\lambda_1,\dots, \lambda_d$ on the diagonal, the modulus $\Omega_n(f,M)_{{p}}$ can be calculated by the following formula (see~\cite{Tim}):
$$
\Omega_n(f,M^{-1})_{{p}}\asymp \sum_{j=1}^d \omega_n^{(j)}(f,\lambda_j^{-1})_p,\quad f\in L_p,\quad 1<p<\infty,
$$
where  $ \omega_n^{(j)}(f,h)_{{p}}=\sup_{|\delta|<h,\, \delta\in \R} \Vert \Delta_{\delta e_j}^n f\Vert_{{p}}$
and $\{e_j\}_{j=1}^d$ denotes the standard basis in $\R^d$.

Let $A$ be a bounded measurable subset of $\R^d$. In what follows, the error of the best approximation of a function $f\in L_{p,w}$ is defined by
$$
E_{A}(f)_{p,w}=\inf\left\{\Vert f-g\Vert_{p,w}\,:\, g\in L_{p,w}\cap L_2,\quad \supp \widehat{g}\subset A\right\}.
$$

\section{Auxiliary results}

An order of approximation by the quasi-projection operators essentially depends on the
compatibility of a function $\vp$ and a distribution/function $\w\vp$.
Assuming that the Fourier transform of $\vp$ and  $\w\vp$ is sufficiently smooth in a neighbourhood  of the origin,
we consider  the following two types of compatibility.

\medskip

\begin{defi}
\label{d1}
 A tempered distribution  $\w\phi$ and a function $\phi$ are said to be
 {\em strictly compatible} if there exists $\delta>0$ such that
 $\overline{\h\phi}(\xi)\h{\w\phi}(\xi)=1$
 a.e. on $\{|\xi|<\delta\}$.
\end{defi}

\begin{defi}
\label{d2+}
 A tempered distribution  $\w\phi$ and a function $\phi$ are said to be
 {\em weakly compatible of order $n\in \N$} if $D^{\beta}(1-\h\phi\overline{\h{\w\phi}})({\bf 0}) = 0$ for all  $[\beta]<n$, $\beta\in\zd_+$.
\end{defi}

\begin{rem}
\emph{ It is well known that the shift-invariant space generated by a function $\phi$ with appropriate decay and smoothness conditions  has approximation order $n$ if and only if
the Strang--Fix conditions of order $n$ are satisfied for $\varphi$, that is $D^\beta\widehat\varphi(l)=0$
whenever $l\in\zd$, $ l\ne\nul$, $[\beta]<n$ (see, e.g.,~\cite{BDR}, \cite{v58}, and~\cite[Ch.~3]{NPS}).
In this paper we restrict ourselves to considering only functions $\varphi$ satisfying the Strang--Fix conditions of arbitrary order $n$.
The condition  $D^{\beta}(1-\h\phi\h{\w\phi})(0) = 0$, $[\beta]<n$, is also a natural requirement for providing approximation order $n$ of quasi-projection operators. This assumption  often appears
(especially in wavelet theory) in other terms, in particular, in terms of polynomial reproducing property (see~\cite[Lemma~3.2]{Jia1}).
It is clear that to provide an infinitely large approximation order, these conditions should be satisfied for any~$n$. Clearly, the latter holds for strictly compatible functions $\vp$ and~$\w\vp$,
 which are used to obtain estimates of approximation via the modulus of continuity of an arbitrary order.
 }
\emph{Supposing that $\h{\w\phi}(\xi)=1$ a.e. on $\{|\xi|<\delta\}$, it is easy to see that the simplest example of $\varphi$ satisfying Definition~\ref{d1} is the tensor product of the $\rm sinc$ functions.}
\end{rem}

\begin{prop}\label{corE}
Let $N\in \z_+$, $\w\phi\in \mathcal{S}_N'$, $\phi\in {L_2}$, $\sum_{k\in\zd}|\h{\phi}(\cdot+k)|^2\in L_\infty$.  Suppose there exists $\delta\in(0,1/2)$ such that
$\h\phi(\xi)=0$ a.e.
 on $\{|l-\xi|<\delta\}$ for all
 $l\in\z^d\setminus\{0\}$, and
$\w\phi$ and $\phi$ be strictly compatible with respect to the parameter $\delta$.
If a function $f\in L_2$ is such that
its Fourier transform is supported in $\{|\xi|<\delta\}$,
then
\be
f=\sum_{k\in\zd} \langle \widehat{f}, \h{\w\phi_{0k}}\rangle \phi_{0k}\quad a.e.
\label{203}
\ee
\end{prop}	

{\bf Proof.}  First of all, we check that the right hand side of~\eqref{203}
belongs to $L_2$.
Set
$$
G(\xi)=\sum_{l\in\zd}\h f(\xi+l)\overline{\h{\widetilde\phi}(\xi+l)}.
$$
By~\cite[Lemma 1]{KS1}, we have that $G\in L_2(\td)$ and $\widehat G(k):=\langle \widehat{f}, \h{\w\phi_{0k}}\rangle$ is its $k$-th Fourier coefficient.
Hence, the series $\sum_{k\in\zd}|\widehat G(k)|^2$ is convergent.
On the other hand, since the system $\{\phi_{0k}\}_{k\in\zd}$ is Bessel (see, e.g., \cite[Remarks~1.1.3 and 1.1.7]{KPS}), we derive that
\begin{equation}\label{Newnew}
   \sum\limits_{k\in\zd}|\langle g,\phi_{0k}\rangle|^2\le B\|g\|_2^2
 \quad\text{for all}\quad g\in L_2,
\end{equation}
where $B=\Vert \sum_{k\in\zd}|\h{\phi}(\cdot+k)|^2\Vert_\infty$.

If now $\Omega$ is a finite subset of $\zd$, by the Riesz representation
theorem, there exists $g\in L_2$ such that $\|g\|_2\le1$ and
$$
\bigg\Vert\sum\limits_{k\in\Omega}
\langle \widehat{f}, \h{\w\phi_{0k}}\rangle
 \phi_{0k}\bigg\Vert_2=
\bigg\Vert\sum\limits_{k\in\Omega}
\widehat G(k) \phi_{0k}\bigg\Vert_2
=\bigg|\sum\limits_{k\in\Omega}\widehat G(k)
\langle\phi_{0k}, g\rangle\bigg|.
$$
Next, using the Cauchy inequality and~\eqref{Newnew}, we get
\ban
\left|\sum\limits_{k\in\Omega}\widehat G(k)
\langle\phi_{0k}, g\rangle\right|\le
\sqrt B\|g\|_2
\lll\sum\limits_{k\in\Omega}|\langle \widehat{f},\widehat{\widetilde\phi_{0k}}
\rangle|^2\rrr^{\frac 12}\le
\sqrt B\lll\sum\limits_{k\in\Omega}|\langle \widehat{f},\widehat{\widetilde\phi_{0k}}
\rangle|^2\rrr^{\frac 12}.
\ean
This implies that the series
$\sum_{k\in\zd} \langle \widehat{f}, \h{\w\phi_{0k}}\rangle \phi_{0k}$  converges as the limit in $L_2$  of the cubic partial sums  and hence its sum belongs to $L_2$.
Using Carleson’s theorem, we have
\begin{equation*}
  \begin{split}
     \mathcal{F}\(\sum_{k\in\zd} \langle \widehat{f}, \h{\w\phi_{0k}}\rangle \phi_{0k}\)(\xi)
&=\sum_{k\in\zd} \widehat G(k)e^{2\pi i(k, \xi)}\h{\phi}(\xi)\\
&=G(\xi)\h{\phi}(\xi)=
\sum\limits_{l\in\,\zd}\h
f(\xi+l)\overline{\h{\w\phi}(\xi+l)}\h{\phi}(\xi)\quad a.e.
   \end{split}
\end{equation*}
The sets $\{|\xi-l|<\delta\}$, $l\in\zd$, are mutually disjoint and
their union contains  the supports of the functions $\h f(\cdot+l)$. If $|\xi-l|<\delta$,
$l\ne0$, then $\h{\phi}(\xi)=0$ and if $|\xi|<\delta$, then $\h{\phi}(\xi)\h{\w\phi}(\xi)=1$. Hence,
$$
\mathcal{F}\(\sum_{k\in\zd} \langle \widehat{f}, \h{\w\phi_{0k}}\rangle \phi_{0k}\)(\xi)=\widehat f(\xi),
$$
 which yields~\eqref{203}~$\Diamond$.
\bigskip

Next we need several basic properties of the modulus of smoothness, which can be proved by  standard arguments using also the inequalities~\eqref{fwa} and~\eqref{fwah}  (see, e.g.,~\cite[Ch.~4]{Nik} and~\cite[Ch.~4]{BeSh}).

\begin{lem}\label{lemmod}
  Let $1\le p\le \infty$, $w\in \mathcal{W}^\alpha$ for some $\alpha>0$, and $n\in \N$. Then for any $f,g\in L_{p,1/w}$ and $\delta\in (0,1)$, we have

\medskip

  \noindent {\rm (i)} $\omega_n(f+g,\delta)_{p,1/w}\le \omega_n(f,\delta)_{p,1/w}+\omega_n(g,\delta)_{p,1/w}$;

\medskip

  \noindent {\rm (ii)} $\omega_n(f,\delta)_{p,1/w}\le C\Vert f\Vert_{p,1/w}$;

\medskip

  \noindent {\rm (iii)} $\omega_n(f,\lambda\delta)_{p,1/w}\le C(1+\lambda)^n\omega_n(f,\delta)_{p,1/w},\, \lambda>0$,

\medskip

\noindent where the constant $C$ does not depend on $f$, $\delta$, and $\lambda>0$.

\end{lem}

 The next result is related to the classical Jackson and Bernstein inequalities (see, e.g.,~\cite[Ch.~7]{BL}); cf. also~Theorem~\ref{jackson}.

\begin{prop}\label{lemKf}
  Let $1\le p\le \infty$, $w\in \mathcal{W}^\alpha$ for some $\alpha>0$, and $n\in \N$. Then for any $f\in L_{p,1/w}$
	 there exists $g\in W_{p,1/w}^n$ such that
  \begin{equation}\label{Kf1}
    \Vert f-g\Vert_{p,1/w}\le C \omega_n(f,1)_{p,1/w}
  \end{equation}
and
  \begin{equation}\label{Kf2}
    \Vert g\Vert_{\dot W_{p,1/w}^n}\le C \omega_n(f,1)_{p,1/w},
  \end{equation}
where the constant $C$ does not depend on $f$ and $w$.
\end{prop}

{\bf Proof.} In the unweighted case, inequalities~\eqref{Kf1} and~\eqref{Kf2}
were proved in~\cite[Theorem 4.12, see also eq. (4.42)]{BeSh} by using the following function $g$:
\begin{equation}\label{stekl1}
  \begin{split}
     g(x)
		=\sum_{k=1}^n (-1)^{k+1}\binom{n}{k} \int_{[0,1]^d}\dots\int_{[0,1]^d} f(x+k (u_1+\dots+u_n))du_1\dots du_n.
   \end{split}
\end{equation}
In the weighted case, inequalities~\eqref{Kf1} and~\eqref{Kf2} can be proved similarly using the same function $g$ and inequalities~\eqref{fwa} and~\eqref{fwah} as well as Lemma~\ref{lemmod} instead of its unweighted counterpart.
~$\Diamond$

\section{Main results}

\subsection{The case of $\vp\in \mathcal{L}_{p,w^*}$ and $\widetilde{\vp}\in \mathcal{L}_{q,w^*}$}

First, we consider the case of strictly compatible functions $\vp$ and $\w\vp$ and give the error estimate in terms of the best approximation.

\begin{theo}\label{lemJack}
Let $1\le p\le\infty$, $w\in \mathcal{W}^\alpha$ for some $\alpha>0$, and $M\in\mathfrak M$.
Suppose
\begin{enumerate}
  \item[1)] $\vp\in \mathcal{L}_{p,w^*}\cap {L}_2$ and $\widetilde{\vp}\in \mathcal{L}_{q,w^*}$,   $1/p+1/q=1$;
  \item[2)] $\vp$ and $\w\vp$
are strictly compatible with respect to the parameter $\delta>0$;
  \item[3)] $\supp \h\vp \subset (-1,1)^d$.
\end{enumerate}
Then for any $f\in L_{p,1/w}$ and $j\in \z_+$
\begin{equation}\label{eqJ1}
    \bigg\|f-\sum_{k\in\zd} \langle f,\w\phi_{jk}\rangle \phi_{jk}\bigg\|_{p,1/w}\le C E_{\{|M^{*-j}\xi|<\delta\}}(f)_{p,1/w},
  \end{equation}
where the constant $C$ does not depend on $f$ and $j$.
\end{theo}

Using Theorem~\ref{lemJack}, the Jackson-type inequality given in Theorem~\ref{jackson}, and Lemma~\ref{lemmod} (iii), we easily obtain the following estimates with moduli of smoothness of arbitrary integer order.

\begin{coro}\label{coro8}
Under the conditions of Theorem~\ref{lemJack}, we have for $p<\infty$ and any $n\in \N$
$$
    \bigg\|f-\sum_{k\in\zd} \langle f,\w\phi_{jk}\rangle \phi_{jk}\bigg\|_{p,1/w}\le C_1 \Omega_n\(f,\sqrt{d}\delta^{-1} M^{-j}\)_{p,1/w}\le C_2\omega_n\(f, \Vert M^{-j}\Vert\)_{p,1/w},
$$
  where $C_1$ and $C_2$ do not depend on $f$ and $j$.
\end{coro}

Now, we consider the case of weakly compatible functions $\vp$ and $\w\vp$. In this case, unlike  Theorem~\ref{lemJack} and Corollary~\ref{coro8}, the error estimate essentially depends on the order of compatibility of $\vp$ and $\w\vp$.

\begin{theo}\label{th2}
Let $1\le p<\infty$, $w\in \mathcal{W}^\alpha$ for some $\alpha>0$,
 $n\in\n$,  and  $M\in\frak M$.
Suppose

\begin{enumerate}
\item[1)] $\phi\in\mathcal{L}_{p, w^*}\cap {L}_2$ and $\w\phi \in\mathcal{L}_{q,w^*}$, $1/p+1/q=1$;

\item[2)] $\h\phi, \h{\w\phi}\in C^{\gamma}(B_\varepsilon)$ for some integer $\gamma>n+d+p\alpha$ and $\varepsilon>0$;

\item[3)] $\w\phi$ and $\phi$ are
 weakly compatible of order $n$;

\item[4)] $\supp \h\vp \subset (-1,1)^d$.
\end{enumerate}

\noindent Then for any $f\in L_{p,1/w}$ and $j\in \z_+$
\begin{equation}
\label{100}
  \bigg\|f-\sum\limits_{k\in\zd}
\langle f,\widetilde\phi_{jk}\rangle \phi_{jk}\bigg\|_{p,1/w}\le C\Omega_n\(f, M^{-j}\)_{p,1/w}\le C \omega_n\(f, \|M^{-j}\|\)_{p,1/w},
\end{equation}
where $C$ does not depend on $f$ and $j$.
\end{theo}

\begin{rem}
\label{rem1}
To determine the approximation order in~(\ref{100}),  one can use a natural approach replacing
$\omega_n\(f,\Vert M^{-j}\Vert\)_{p,1/w}$ by $\omega_n\(f,\Vert M^{-1}\Vert^{j}\)_{p,1/w}$. However this is not
good because there exist matrices $M\in\frak M$ such that $\|M^{-1}\|\ge1$.
A better way is to use~\eqref{00++}, which yields
$$
  \bigg\|f-\sum\limits_{k\in\zd}
\langle f,\widetilde\phi_{jk}\rangle \phi_{jk}\bigg\|_{p,1/w}\le C\omega_n\(f, \vartheta^{-j}\)_{p,1/w},
$$
where  $\vartheta$  is a positive number smaller (in absolute value)
than any eigenvalues of $M$.  If $M$ is
an isotropic matrix and $\lambda$ is one of its eigenvalues (e.g., $M=\lambda I_d$), then one can take $\vartheta=|\lambda|$.
If $M$ is a diagonal matrix and $\lambda$   is the smallest (in absolute value) diagonal element, then one can take $\vartheta=|\lambda|$.
\end{rem}

The proof of Theorem~\ref{th2} (see below) provides the following Jackson-type theorem in the weighted spaces $L_{p,1/w}$.

\begin{theo} \emph{(Jackson inequality)}
\label{jackson}
Let  $1\le p<\infty$, $w\in{\cal W}^\alpha$ for some $\alpha>0$, $n\in\n$, and $M\in \mathfrak{M}$. Then for any  $f\in L_{p,1/w}$ and $j\in \Z_+$
$$
E_{\{|M^{*-j}\xi|<\sqrt{d}\}}(f)_{p,1/w}\le C\Omega_n\(f, M^{-j}\)_{p,1/w} \le C\omega_n\(f,\Vert M^{-j}\Vert\)_{p,1/w},
$$
where $C$ does not depend on $f$ and $j$.
\end{theo}

To prove the above theorems, we will use the following statements.
	
\begin{prop}\label{uns2}
  Let $1\le p\le \infty$ and $w\in \mathcal{W}^\alpha$ for some $\alpha>0$.
 If $\vp\in \mathcal{L}_{p,w^*}$ and $c=\{c_k\}_{k\in \Z^d}\in \ell_{p,w}$, then
$$
\bigg\|\sum_{k\in\zd} c_k \phi_{0k}\bigg\|_{p,w}\le  \Vert\vp\Vert_{\mathcal{L}_{p,w^*}}\|c\|_{\ell_{p,w}}.
$$
\end{prop}

\begin{prop}\label{uns2'}
  Let $1\le p\le \infty$, $w\in \mathcal{W}^\alpha$ for some $\alpha>0$,
  $\vp\in \mathcal{L}_{p,w^*}$, and $c=\{c_k\}_{k\in \Z^d}\in \ell_{p,1/w}$. Then
$$
\bigg\|\sum_{k\in\zd} c_k \phi_{0k}\bigg\|_{p,1/w}\le  \Vert\vp\Vert_{\mathcal{L}_{p,w^*}}\|c\|_{\ell_{p,1/w}}.
$$
\end{prop}

\begin{prop}\label{prop1Lq}
Let $1\le p\le \infty$, $w\in \mathcal{W}^\alpha$ for some $\alpha>0$, $f\in L_{p,1/w}$, and $\vp\in \mathcal{L}_{q,w^*}$, $1/p+1/q=1$.
Then
$$
  \left(\sum_{k\in\zd} \bigg|\frac{\langle f,{{\phi}_{0k}}\rangle}{w(k)}\bigg|^p\right)^\frac 1p\le \Vert \vp\Vert_{\mathcal{L}_{q,w^*}}\|f\|_{p,1/w}.
$$
\end{prop}

The above three propositions can be proved repeating step-by-step the proofs of Propositions~2, 4, and~5 in~\cite{NU}
respectively.~$\Diamond$

\begin{coro}\label{lemBound}
  Let $1\le p\le\infty$ and $w\in \mathcal{W}^\alpha$ for some $\alpha>0$. If $\vp\in \mathcal{L}_{p,w^*}$ and $\widetilde{\vp}\in \mathcal{L}_{q,w^*}$, $1/p+1/q=1$, then for  any $f\in L_{p,1/w}$
  \begin{equation}\label{eq0}
    \bigg\|\sum_{k\in\zd} \langle f,\w\phi_{0k}\rangle \phi_{0k}\bigg\|_{p,1/w}
		\le  \Vert\vp\Vert_{\mathcal{L}_{p,w^*}}\Vert\w\vp\Vert_{\mathcal{L}_{q,w^*}} \Vert f\Vert_{p,1/w}.
  \end{equation}
 \end{coro}

{\bf Proof.} The proof of~\eqref{eq0} directly follows from   Propositions~\ref{uns2'} and~\ref{prop1Lq}.~$\Diamond$

\bigskip

{\bf Proof of Theorem~\ref{lemJack}.}
Without loss of generality, we can assume that $\delta$  is sufficiently small such that $\h\vp(\xi)= 0$ a.e. on $|\xi-l|<\delta$ for all $l\in\zd\setminus\{\nul\}$.

If $g\in L_{p,1/w}\cap L_2$ and $\supp \widehat{g}\subset \{|\xi|<\delta\}$, then,
due to Proposition~\ref{corE}, we have
\begin{equation}\label{NNN}
  g=\sum_{k\in\zd} \langle g,\w\phi_{0k}\rangle \phi_{0k}.
\end{equation}
Corollary~\ref{lemBound} and \eqref{NNN} imply that
\begin{equation}\label{eqJ2}
  \begin{split}
    \bigg\|f-\sum_{k\in\zd} \langle f,\w\phi_{0k}\rangle \phi_{0k}\bigg\|_{p,1/w}
    &\le \Vert f-g\Vert_{p,1/w}+\bigg\|\sum_{k\in\zd} \langle g-f,\w\phi_{0k}\rangle \phi_{0k}\bigg\|_{p,1/w}\\
    &\le (1+ \Vert\vp\Vert_{\mathcal{L}_{p,w^*}}\Vert\w\vp\Vert_{\mathcal{L}_{q,w^*}})
		\Vert f-g\Vert_{p,1/w}.
  \end{split}
\end{equation}

Let  now $j\in\z_+$ be fixed, $G$ be a function in $L_{p,1/w}\cap L_2$
such that $\supp \widehat G\subset \{|{M^*}^{-j}\xi|<\delta\}$
 and
\begin{equation}\label{eq*}
  \Vert f-G\Vert_{p, 1/w}\le 2 E_{\{|{M^*}^{-j}\xi|<\delta\}}(f)_{p,1/w}.
	\end{equation}
	Set $g(x)=G(M^{-j}x)$. Obviously, $g\in L_{p,1/w(M^{-j}\cdot)}\cap L_2$ and $\supp \widehat{g}\subset \{|\xi|<\delta\}$.
 Thus, after the change of variable,  using~\eqref{eqJ2} with $f(M^{-j}\cdot)$ instead of $f$ and $w(M^{-j}\cdot)$ instead of $w$, we have
\ban
\bigg\|f-\sum_{k\in\zd} \langle f,\w\phi_{jk}\rangle \phi_{jk}\bigg\|_{p,1/w}
    =m^{-j/p}\bigg\|f(M^{-j}\cdot)-\sum_{k\in\zd} \langle f(M^{-j}\cdot),\w\phi_{0k}\rangle \phi_{0k}\bigg\|_{p,1/w(M^{-j}\cdot)}
 \\
	\le\big(1+ \Vert\vp\Vert_{\mathcal{L}_{p,w^*(M^{-j}\cdot)}}\Vert\w\vp\Vert_{\mathcal{L}_{q,w^*(M^{-j}\cdot)}}\big)\,
m^{-j/p}	\Vert f(M^{-j}\cdot)- g\Vert_{p,1/w(M^{-j}\cdot)}.
  \ean
To prove~\eqref{eqJ1}, it remains to note that
	$$
		m^{-j/p}\Vert f(M^{-j}\cdot)-g\Vert_{p,1/w(M^{-j}\cdot)}=
			\Vert f-G\Vert_{p,1/w},
	$$
apply~\eqref{eq*}, and take into account that
$$
\Vert\vp\Vert_{\mathcal{L}_{p,w^*(M^{-j}\cdot)}}\le C'\Vert\vp\Vert_{\mathcal{L}_{p,w^*}}
\quad
\text{and}
\quad
\Vert\w\vp\Vert_{\mathcal{L}_{q,w^*(M^{-j}\cdot)}}\le C'\Vert\w\vp\Vert_{\mathcal{L}_{q,w^*}},
$$
which proves the theorem.~$\Diamond$

\medskip

To prove  Theorem~\ref{th2} we need the following lemmas.

\begin{lem}\label{lem1}
Let $1\le p\le \infty$, $w\in \mathcal{W}^\alpha$ for some $\alpha>0$, and
 $n\in\n$.
Suppose
\begin{enumerate}
  \item[1)]  $\psi\in\mathcal{L}_{p, w^*}$;
  \item[2)]  $\h{\w\psi}\in C^{\gamma}(\rd)$ for some $\gamma>n+d+p\alpha$ and $\h{\w\psi}$ is compactly supported;
  \item[3)]  $D^{\beta}\h{\w\psi}(\nul) = 0$ for all $[\beta]<n$, $\beta\in \zd_+$.
\end{enumerate}
Then for any $f\in W_{p,1/w}^n$
\begin{equation*}
   \bigg\|\sum\limits_{k\in\zd}
\langle f,\widetilde\psi_{0k}\rangle \psi_{0k}\bigg\|_{p, 1/w}\le  c_w^2C\|\psi\|_{\mathcal{L}_{p, w^*}} \|f\|_{\dot W^n_{{p,1/w}}},
\end{equation*}
where $C$ does not depend on $f$.
\end{lem}

\textbf{Proof.}
Let $k\in\zd$ and $y\in [-1/2,1/2]^d-k$. Since $D^{\beta}\h{\w\psi}(\nul) = 0$ whenever $[\beta]<n$, $\beta\in \Z_+^d$, we have
$$
\int\limits_{\rd}x^\beta\widetilde\psi_{0k}(x)\,dx=0,\quad [\beta]< n,\quad
 \beta\in\zd_+.
$$
Hence, due to Taylor's formula with the integral remainder,
\begin{equation*}
  \begin{split}
      |\langle f,\w\psi_{0k}\rangle |&=
    \bigg|\int\limits_{\rd}f(x)\overline{\w\psi_{0k}(x)}\,dx\bigg|
    \\
    &=\bigg|\int\limits_{\rd}\, \overline{\w\psi_{0k}(x)}\bigg(\sum\limits_{\nu=0}^{n-1}\frac{1}{\nu!}\lll (x_1-y_1)\partial_{1}
    + \dots + (x_d-y_d)\partial_{d} \rrr^\nu f(y)
    \\
    &\quad\quad\quad\quad+\int\limits_0^1\frac{(1-t)^{n-1}}{(n-1)!} \Big((x_1-y_1)\partial_{1}
		+ \dots + (x_d-y_d)\partial_{d} \Big)^n f(y+t(x-y))\,dt\bigg)\,dx \bigg|
    \\
    &\le \int\limits_{\rd}|x-y|^n|\w\psi_{0k}(x)|\,\int\limits_0^1\sum\limits_{[\beta]=n}|D^{\beta}f(y+t(x-y))|\,dt\,dx\,.
   \end{split}
\end{equation*}
From this, using  H\"older's inequality and taking into account that		
$$
|\w\psi_{0k}(x)|\le\frac{C_1 }{{(1+|x+k|)^\gamma}}\le \frac{C_2 }{(1+|x-y|)^\gamma},
$$
we obtain
\begin{equation}\label{dlin}
  \begin{split}
          |\langle f,\w\psi_{0k}\rangle |&\le
   C_2   \int\limits_{\rd}\frac{| x-y|^n}{(1+|x-y|)^\gamma}
    \int\limits_0^1\sum\limits_{[\beta]=n}|D^{\beta}f(y+t(x-y))|\,dt \,dx
    \\
    &\le C_2
\bigg(\int\limits_{\rd}\frac{| x-y|^{n}}{(1+|x-y|)^\gamma}\,dx\bigg)^{1/q}
\\
    &\quad\quad\quad\quad\times\Bigg(\int\limits_{\rd}\,\frac{|x-y|^{n}}{(1+|x-y|)^\gamma}\bigg(\int\limits_0^1\,
\sum\limits_{[\beta]=n}|D^{\beta}f(y+t(x-y))|\,dt\bigg)^p dx\Bigg)^{1/p}
\\
	&=C_3
    \Bigg(\int\limits_{\rd}\,\frac{|u|^n}{(1+|u|)^\gamma}\int\limits_0^1\,\sum\limits_{[\beta]=
	n}|D^{\beta}f(y+tu)|^p\,dt \,du\Bigg)^{1/p}.
   \end{split}
\end{equation}
It follows from~\eqref{dlin}, properties of $w$, and  Proposition~\ref{uns2'} that
\begin{equation*}
\begin{split}	
        & \bigg\|\sum\limits_{k\in\,\zd}\langle
    f,\w\psi_{0k}\rangle\psi_{0k}\bigg\|_{p,1/w}^p
    \le
    \|\psi\|^p_{\mathcal{L}_{p, w^*}}\sum_{k\in\zd}\left|\frac{\langle f,\w\psi_{0k}\rangle}{w(k)} \right|^p
			=\|\psi\|_{\mathcal{L}_{p, w^*}}^p\sum_{k\in\zd}\int\limits_{[-\frac12,\frac12]^d-k}\,dy
		\left|\frac{\langle f,\w\psi_{0k}\rangle}{w(k)} \right|^p
\\
		&\le C_4 \|\psi\|^p_{\mathcal{L}_{p, w^*}}
		\sum_{k\in\zd}\int\limits_{[-\frac12,\frac12]^d-k}
    \int\limits_{\rd}\,\frac{|u|^n}{(1+|u|)^\gamma}\int\limits_0^1
		\left|\frac{w(y)w^*(tu)}{w(k)}\right|^p
		\sum\limits_{[\beta]=n}\left|\frac{D^{\beta}f(y+tu)}{w(y+tu)}\right|^p\,dt\,du\,dy
\\
    &\le C_4 c_w^p\|\psi\|^p_{\mathcal{L}_{p, w^*}}
		\sum_{k\in\zd}\int\limits_{[-\frac12,\frac12]^d-k}(w^*(y+k))^p\times
\\
			&\qquad\qquad\qquad\qquad\qquad\qquad\int\limits_{\rd}\,\frac{|u|^n }{(1+|u|)^\gamma}
		 \int\limits_0^1\,		(1+|tu|^2)^{{p\alpha}/2}
  	\sum\limits_{[\beta]=n}\left|\frac{D^{\beta}f(y+tu)}{w(y+tu)}\right|^p\,dt\,du\,dy\\
\end{split}
\end{equation*}
Next, using Fubini's Theorem, the fact that
$(1 + |tu|^2)^{p\alpha/2} \le (1 + |u|^2)^{p\alpha/2}$, we derive from the above estimate that
\begin{equation*}
\begin{split}
          \bigg\|\sum\limits_{k\in\,\zd}\langle
    f,\w\psi_{0k}\rangle\psi_{0k}\bigg\|_{p,1/w}^p
     &\le C_5 c_w^{2p}\|\psi\|^p_{\mathcal{L}_{p, w^*}}
    \int\limits_{\rd} \frac{|u|^n (1+|u|^2)^{p\alpha/2}du}{(1+|u|)^\gamma}
		\int\limits_{\rd}\sum\limits_{[\beta]=n}\left|\frac{D^{\beta}f(y)}{w(y)}\right|^p\,dy
			\\
    &\le C_6 c_w^{2p}\|\psi\|^p_{\mathcal{L}_{p, w^*}}
		\int\limits_{\rd}\,\frac{du}{(1+|u|)^{\gamma-n-p\alpha}}\|f\|_{\dot W^n_{p,1/w}}^p
		\le C_7 c_w^{2p} \|\psi\|^p_{\mathcal{L}_{p, w^*}} \|f\|_{\dot W^n_{{p,1/w}}}^p.
\end{split}
\end{equation*}
This proves the lemma.~$\Diamond$

\medskip

\begin{lem}\label{lem2}
Let $w$, $n$, $\gamma$, $\phi$, and $\w\phi$ be as in Theorem~\ref{th2}.
Then for any $f\in W_{p,1/w}^n$, we have
\begin{equation}\label{101}
   \bigg\|f-\sum\limits_{k\in\zd}
\langle f,\widetilde\phi_{0k}\rangle \phi_{0k}\bigg\|_{p, 1/w}\le  \Upsilon(w^*) \|f\|_{\dot W^n_{{p,1/w}}},
\end{equation}
where the functional $\Upsilon$ is independent of $f$ and
$
 \Upsilon(w^*(M^{-j}\cdot))\le C\Upsilon(w^*)
$
for all $j\in\z_+$ and the constant $C$ depends only on
$d$, $p$, $\phi$, $\w \phi$, $\alpha$, $C'$, and $c_w$.
\end{lem}

\textbf{Proof.}
Due to Corollary~\ref{lemBound}, without loss of generality, we may assume that $f\in L_2$.
Choose $0<\delta'<\delta''<1/2$ such that $\h\phi(\xi)\ne0$ on $\{|\xi|\le\delta'\}$.
Set
$$
F(\xi)=
 \begin{cases}
	 \displaystyle\frac{1-\overline{\h\phi(\xi)}\h{\w\phi}(\xi)}{\overline{\h\phi(\xi)}}  &\mbox{if $|\xi|\le\delta'$,}
		 \\
	\displaystyle 0
	 &\mbox{if $|\xi|\ge\delta''$}
	 		\end{cases}
	$$
and extend this function  such that $F\in C^{\gamma}(\rd)$. Define $\w\psi$ by
$\h{\w\psi}=F$. Obviously, the function $\w\psi$ is continuous
and $\w\psi(x)=O(|x|^{-\gamma})$ as $|x|\to \infty$, where
$\gamma>\alpha+d$, which yields that $\w\psi\in \mathcal{L}_{\infty, w^*}$, a fortiori $\w\psi\in \mathcal{L}_{q, w^*}$.
On the other hand,  all assumptions of Lemma~\ref{lem1}
with $\psi=\phi \in \mathcal{L}_{p,w^*}$ are satisfied. Hence,
\be
\label{108}
\bigg\Vert\sum_{k\in\zd}\langle f, \w\psi_{0k}\rangle \vp_{0k}\bigg\Vert_{p,1/w}\le
C_0  c_w^2 \|\phi\|_{\mathcal{L}_{p, w^*}}
\|f\|_{\dot W_{p,1/w}^n},
\ee
where $C_0$ depends on $d$, $p$, $\w\psi$, and $\alpha$. Using~\eqref{108}, we derive
\begin{equation}\label{108+++}
  \bigg\|f-\sum\limits_{k\in\zd}
\langle f,\widetilde\phi_{0k}\rangle \phi_{0k}\bigg\|_{p, 1/w}\le C_0  c_w^2 \|\phi\|_{\mathcal{L}_{p, w^*}}
\|f\|_{\dot W_{p,1/w}^n}+\bigg\Vert f-\sum_{k\in\zd}\langle f, \w\phi_{0k}+\w\psi_{0k}\rangle \vp_{0k}\bigg\Vert_{p,1/w}.
\end{equation}

It remains to estimate the second summand in the right-hand side of~\eqref{108+++}.
For this, we will use the Meyer wavelets. Let $\theta$ be the Meyer scaling function
 (see, e.g.~\cite[Sec.~1.4] {NPS}). This function is band-limited, its Fourier transform  is
infinitely differentiable, supported in $[-2/3,2/3]$, and equals 1 on the interval $[-1/3, 1/3]$;
the integer translates of $\theta$ form an orthonormal system. Set
$$
\Phi(x)=\prod_{l=1}^d\theta(x_l),\quad x\in\rd.
$$
It is well known  (see, e.g.,~\cite[Sec.~2.1] {NPS}) that $\Phi$ generates a separable MRA in $L_2$
with respect to the  matrix $2I_d$ and  the corresponding
wavelet functions $\Psi^{(\nu)}$, $\nu=1,\dots, 2^d-1,$ such that
for every $j\in \z$ the functions $\Phi_{jk}=2^{jd/2}\Phi(2^j\cdot+k)$ and $\Psi^{(\nu)}_{ik}= 2^{id/2}\Phi^{(\nu)}(2^i\cdot+k)$,
$k\in\zd$, $i\ge j$, $\nu=1,\dots, 2^j-1$, form an orthonormal basis for $L_2$.
It follows that
\be
\label{105}
f=\sum_{k\in\zd}\langle f, \Phi_{jk}\rangle \Phi_{jk} +\sum_{\nu=1}^{2^d-1}\sum_{i=j}^{\infty}
\sum_{k\in\zd}\langle f, \Psi^{(\nu)}_{ik}\rangle \Psi^{(\nu)}_{ik}.
\ee

On the other hand, the functions $\h\Psi^{(\nu)}$ are infinitely differentiable and compactly supported. It follows that
(see~\cite[Theorem 1.7.7 and Sec.~1.4] {NPS})
\be
\label{103}
\int\limits_{\rd}y^\beta\Psi^{(\nu)}(y)\,dy=0 \quad\text{for all}\quad  \beta\in\zd_+.
\ee
Thus, all assumptions of  Lemma~\ref{lem1} with $\psi=\w\psi=\Psi^{(\nu)}$
 are satisfied. Hence
$$
\bigg\Vert\sum_{k\in\zd}\langle f, \Psi^{(\nu)}_{0k}\rangle \Psi^{(\nu)}_{0k}\bigg\Vert_{p,1/w}\le
C_1  c_w^2 \|\Psi\|_{\mathcal{L}_{p, w^*}}
\|f\|_{\dot W_{p,1/w}^n},\quad \nu= 1,\dots, 2^d-1,
$$
where $C_1$ depends on $d$, $p$, $\Psi$ and $\alpha$.   For every $i\in\z$ and  $\nu= 1,\dots, 2^d-1$,
after the change of variable, taking into account that $w(2^{-i}\cdot)\in \mathcal{W}^\alpha$ and
$$
w^*(2^{-i}x)\le c_w \max\{1, 2^{-i\alpha}\}(1+|x|^2)^{\alpha/2},
$$
  we have
\be
\label{114}
\bigg\Vert\sum_{k\in\zd}\langle f, \Psi^{(\nu)}_{ik}\rangle \Psi^{(\nu)}_{ik}\bigg\Vert_{p,1/w}\le
C_1 c_w^2 \max\{1, 2^{-2i\alpha}\}\|\Psi\|_{\mathcal{L}_{p, w^*(2^{-i}\cdot)}}
2^{-id/p}\|f(2^{-i}\cdot)\|_{\dot W_{p,1/w(2^{-i}\cdot)}^n}.
\ee

 Choose $j\in\z$ such that $2^{j}<\delta'/\sqrt d$ and set
 $G=\sum_{k\in\zd}\langle f, \Phi_{j k}\rangle \Phi_{j k}$. Note that $j=j(\phi, \w\phi)$.
Since $\supp \h\Phi\subset [-1,1]^d$,
we have
\be
\supp\h G\subset2^{j}[-1,1]^d\subset \delta' B_1.
\label{106}
\ee

By construction,  $\overline{\h\phi(\xi)}(\h{\w\phi}(\xi)+\h{\w\psi}(\xi))=1$ whenever $|\xi|\le\delta'$.
It follows from~(\ref{106}) and Proposition~\ref{corE} that
\be
\label{107}
G=\sum_{k\in\zd}\langle G, \w\phi_{0k}+\w\psi_{0k}\rangle \vp_{0k}.
\ee
Since $\w\phi, \w\psi\in \mathcal{L}_{q, w^*}$, $\phi\in \mathcal{L}_{p,w^*}$, and $f, G\in L_{p,1/w}$,
due to Corollary~\ref{lemBound}, we derive
\be
\label{108+}
     \bigg\Vert\sum_{k\in\zd}\langle f-G, \w\phi_{0k}+\w\psi_{0k}\rangle \vp_{0k}\bigg\Vert_{p,1/w}\le
 \Vert\phi\Vert_{\mathcal{L}_{p,w^*}}\Vert\w \phi+\w \psi\Vert_{\mathcal{L}_{q,w^*}}\|f-G\|_{p,1/w}.
\ee
Combining  (\ref{108+}), (\ref{105}), and (\ref{107}), we obtain
\begin{equation}
	   \begin{split}
		\label{102}
&\bigg\|f-\sum_{k\in\zd}\langle f, \w\phi_{0k}+\w\psi_{0k}\rangle \vp_{0k}\bigg\|_{p,1/w}
=\bigg\|f-G-\sum_{k\in\zd}\langle f-G, \w\phi_{0k}+\w\psi_{0k}\rangle \vp_{0k}\bigg\|_{p,1/w}
\\
&\le \bigg(1+ \Vert\phi\Vert_{\mathcal{L}_{p,w^*}}\Vert\w \phi+\w \psi\Vert_{\mathcal{L}_{q,w^*}}\bigg) \Vert f-G\Vert_{p,1/w}
\\
&=\bigg(1+ \Vert\phi\Vert_{\mathcal{L}_{p,w^*}}\Vert\w \phi+\w \psi\Vert_{\mathcal{L}_{q,w^*}}\bigg)
\bigg\|\sum_{\nu=1}^{2^d-1}\sum_{i=j}^{\infty}
\sum_{k\in\zd}\langle f, \Psi^{(\nu)}_{ik}\rangle \Psi^{(\nu)}_{ik}\bigg\|_{p,1/w}.
	   \end{split}
\end{equation}
It follows from~(\ref{114}) that
\begin{equation}
	   \begin{split}
		\label{104}
&\bigg\|\sum_{\nu=1}^{2^d-1}\sum_{i=j}^{\infty}
\sum_{k\in\zd}\langle f, \Psi^{(\nu)}_{ik}\rangle \Psi^{(\nu)}_{ik}\bigg\|_{p,1/w}
\\
&\le C_1 c_w^2 \max\{1, 2^{-2j\alpha}\} \sup\limits_{i\ge j}\|\Psi\|_{\mathcal{L}_{p, w^*(2^{-i}\cdot)}}
 \sum_{\nu=1}^{2^d-1}\sum_{i=j}^{\infty}2^{-i d/p} \|f(2^{-i}\cdot)\|_{\dot W_{p,1/w((2^{-i}\cdot)}^n}
\\
&= C_1 c_w^2 \max\{1, 2^{-2j\alpha}\} \sup\limits_{i\ge j}\|\Psi\|_{\mathcal{L}_{p, w^*(2^{-i}\cdot)}}
\lll\sum_{\nu=1}^{2^d-1}\sum_{i=j}^{\infty}2^{-ni}\rrr\|f\|_{\dot W_{p,1/w}^n}.
	   \end{split}
\end{equation}
Since $w^*(x/2)\le c_w(1+|x|^2)^{\alpha/2}$ and $\h\Psi$ is infinitely differentiable, we have
$$
\sup\limits_{i\ge j}\|\Psi\|_{\mathcal{L}_{p, w^*(2^{-i}\cdot)}} \le C_2.
$$
Combining  this with (\ref{102}), (\ref{104}), and (\ref{108+++}), we complete the proof of~(\ref{101})
with
$$
\Upsilon(w^*)=C_3\max\bigg\{1, \Vert\phi\Vert_{\mathcal{L}_{p,w^*}}\Vert\w \phi+\w \psi\Vert_{\mathcal{L}_{q,w^*}},
\Vert\phi\Vert_{\mathcal{L}_{p,w^*}}\bigg\},
$$
where $C_3$ depends only on
$d$, $p$, $\phi$, $\w \phi$, $\alpha$, and $c_w$.~$\Diamond$

\bigskip

\textbf{Proof of Theorem~\ref{th2}.}
Let $g\in W_{p,1/w}^n$ be defined by~\eqref{stekl1}.
Using Corollary~\ref{lemBound}, Lemma~\ref{lem2}, and Proposition~\ref{lemKf}, we derive
\begin{equation}
\label{120}
  \begin{split}
       \bigg\|f-&\sum\limits_{k\in\zd}
\langle f,\widetilde\phi_{0k}\rangle \phi_{0k}\bigg\|_{p,1/w}\\
&\le \Vert f-g\Vert_{p,1/w}+\bigg\|g-\sum\limits_{k\in\zd}
\langle g,\widetilde\phi_{0k}\rangle \phi_{0k}\bigg\|_{p,1/w}+\bigg\|\sum\limits_{k\in\zd}
\langle f-g,\widetilde\phi_{0k}\rangle \phi_{0k}\bigg\|_{p,1/w}\\
&\le C_1 \Upsilon(w^*)\(\Vert f-g\Vert_{p,1/w}+\Vert g\Vert_{\dot W_{p,1/w}^n}\)
\le C_2 \Upsilon(w^*) \omega_n\(f,1\)_{p,1/w},
   \end{split}
\end{equation}
where $C_2$ depends only on $d$, $p$, $\phi$, $\w \phi$, $\alpha$, and $c_w$.
This yields~\eqref{100} for $j=0$.

Consider now an arbitrary $j\in\z_+$. After the change of variables, we have
$$
 \bigg\|f-\sum\limits_{k\in\zd}
\langle f,\widetilde\phi_{jk}\rangle \phi_{jk}\bigg\|_{p,1/w}=m^{-j/p} \bigg\|f(M^{-j}\cdot)-\sum\limits_{k\in\zd}
\langle f(M^{-j}\cdot),\widetilde\phi_{0k}\rangle \phi_{0k}\bigg\|_{p,1/w(M^{-j}\cdot)}.
$$
Since $w(M^{-j}\cdot)\in {\cal W}^{\alpha}$ and $\Upsilon\((w^*(M^{-j}\cdot)\)\le C_3\Upsilon(w^*)$,
it follows from~(\ref{120}) that
\begin{equation*}
  \begin{split}
  \bigg\|f-\sum\limits_{k\in\zd} \langle f,\widetilde\phi_{jk}\rangle \phi_{jk}\bigg\|_{p,1/w}
  &\le C_4 \Upsilon(w^*) m^{-j/p}\omega_n\(f(M^{-j}\cdot),1\)_{p,1/w(M^{-j}\cdot))}\\
  &\le C_4 \Upsilon(w^*) \Omega_n\(f, M^{-j}\)_{p,1/w}\le C_5\omega_n\(f, \|M^{-j}\|\)_{p,1/w},
   \end{split}
\end{equation*}
which proves the theorem.~$\Diamond$

\medskip

\textbf{Proof of Theorem~\ref{jackson}.}
Choose  a function $g\in \mathcal{S}$ such that
\be
\label{111}
\|f-g\|_{p,1/w}\le \Omega_n\(f,M^{-j}\)_{p,1/w}
\ee
and set
$$
P=\sum_{k\in\zd}\langle g, {\Phi}_{jk}\rangle {\Phi}_{jk},
$$
where $\Phi$ is the scaling function of the separable
Meyer MRA (see the proof of Lemma~\ref{lem2}). All conditions of Theorem~\ref{th2} with $\phi=\w\phi={\Phi}$ are satisfied, and hence
$$
\|g-P\|_{p,1/w}\le C_1\Omega_n\(g, M^{-j}\)_{p,1/w}
$$
 Thus, taking into account that
$P\in L_{p,1/w}\cap L_2$ and $\supp \widehat{P}\subset {\{|M^{*-j}\xi|<\sqrt{d}\}}$, we obtain
$$
E_{\{|M^{*-j}\xi|<\sqrt{d}\}}(f)_{p,1/w}\le\|f-P\|_{p,1/w}
\le C_1\Omega_n\(f, M^{-j}\)_{p,1/w}+C_2\|f-g\|_{p,1/w}.
$$
It remains to combine this with~(\ref{111}).~$\Diamond$

\subsection{The case of $\vp\in \mathcal{B}$ and  $\w\phi \in\mathcal{L}_{q,w^*}$}
\label{S5.2}

Let $\mathcal{X}$ be a collection of  bounded sets in $\R^d$ and let $w$ be a nonnegative, locally integrable function. We say that $w$ belongs to $A_p(\R^d, \mathcal{X})$ for some $1<p<\infty$ if there is a constant $c$ such that
$$
\(\frac1{\mes I} \int_I w(x)dx\)\(\frac1{\mes I}\int_I w(x)^{-1/(p-1)} dx\)^{p-1}\le c
$$
for any $I\in \mathcal{X}$.

Now let $\mathcal{Q}_d$ and $\mathcal{R}_d$ denote the collection of all $d$-dimensional cubes and all $d$-dimensional rectangles with sides parallel to the coordinate axes, correspondingly. Then $A_p(\R^d, \mathcal{Q}_d)$ is the classical Muckenhoupt class $A_p(\R^d)$. If $d=1$, then $A_p(\R^1)=A_p(\R^1, \mathcal{R}_1)$. At the same time, for $d>1$, we have $A_p(\R^d,\mathcal{R}_d)\subsetneqq A_p(\R^d)$. Recall also the fact that $|x|^\alpha \in A_p(\R^d)$ for $-d<\alpha<d(p-1)$ while $|x|^\alpha\in A_p(\R^d, \mathcal{R}_d)$ for $-1<\alpha<p-1$ (see, e.g.,~\cite{Kurtz}).

In what follows, for simplicity we denote ${\mathcal{A}}_p={\mathcal{A}}_p(\rd)=A_p(\R^d, \mathcal{R}_d)$.

In the results formulated in the previous sections, we suppose that the weight $w$ belongs to the class ${\cal W}^{\alpha}$. In the next results, we will in addition suppose that $w^{-p}\in {\mathcal{A}}_p$ (or, equivalently, $w^q\in {\mathcal{A}}_q$ with $1/q+1/p=1$). A model example of such a weight is
$$
w_\alpha(x)=(1+|x|^2)^{\alpha/2},\    \ 0<\alpha<1/p,\  \ 1<p<\infty .
$$

\begin{theo}\label{lemJack+++}
Let $1< p<\infty$, $w\in {\cal W}^{\alpha}$ for some $\alpha \in (0,1)$, $w^{-p}\in {\mathcal{A}}_p$, and let  $M \in \mathfrak{M}$ be a diagonal  matrix.
Suppose
\begin{enumerate}
  \item[1)] $\vp\in \mathcal{B}$ and $\widetilde{\vp}\in \mathcal{L}_{q,w^*}$,   $1/p+1/q=1$;
  \item[2)] $\vp$ and $\w\vp$
are strictly compatible with respect to the parameter $\delta>0$;
  \item[3)]  $\supp \h\vp \subset (-1,1)^d$.
\end{enumerate}
Then for any $f\in L_{p,1/w}$ and $j\in \z_+$
\begin{equation*}
    \bigg\|f-\sum_{k\in\zd} \langle f,\w\phi_{jk}\rangle \phi_{jk}\bigg\|_{p,1/w}\le C E_{\{|M^{-j}\xi|<\delta\}}(f)_{p,1/w},
  \end{equation*}
where the constant $C$ does not depend on $f$ and $j$.
\end{theo}

In the case of weakly compatible functions $\vp$ and $\w\vp$, we have the following result.


\begin{theo}\label{th4}
 Let $1< p<\infty$, $w\in {\cal W}^{\alpha}$ for some $\alpha \in (0,1)$, $w^{-p}\in {\mathcal{A}}_p$, $n\in\n$, and let  $M \in \mathfrak{M}$ be a diagonal  matrix.
Suppose

\begin{enumerate}
  \item[1)]  $\phi\in\mathcal{B}$ and $\w\phi \in\mathcal{L}_{q,w^*}$, $1/p+1/q=1$;
  \item[2)]  $\h\phi, \h{\w\phi}\in C^{\gamma}(B_\varepsilon)$ for some integer $\gamma>n+d+p\alpha$ and $\varepsilon>0$;
  \item[3)]  $\w\phi$ and $\phi$ are weakly compatible of order $n$;
   \item[4)]  $\supp \h\vp \subset (-1,1)^d$.
\end{enumerate}

\noindent Then for any $f\in L_{p,1/w}$ and $j\in \Z_+$, we have
\begin{equation*}
  \bigg\|f-\sum\limits_{k\in\zd}
\langle f,\widetilde\phi_{jk}\rangle \phi_{jk}\bigg\|_{p,1/w}\le
C \Omega_n\(f,M^{-j}\)_{p,1/w}
\le
C \omega_n\(f, |\lambda|^{-j}\)_{p,1/w},
\end{equation*}
where  $\lambda$ is the smallest (in absolute value) diagonal element of $M$ and $C$ does not depend on $f$  and~$j$.
\end{theo}

To prove Theorems~\ref{lemJack+++} and~\ref{th4}, we will use the following auxiliary results.

\begin{lem}\label{lemAp} {\rm (See~\cite[p.~453--454]{GCRF} and \cite{Kurtz}.)}
Let $1<p<\infty$ and $d\ge 2$. The following assertions are equivalent:

\begin{enumerate}
\item[1)] $w\in {\mathcal{A}}_p(\rd)$;

\item[2)] There exists a constant $C$ such that for almost every fixed vector  $(x_1,\dots,x_{j-1},x_{j+1},\dots,x_d)\in \R^{d-1}$ and any interval $I\subset \R^1$ one has
$$
\(\frac1{\mes I} \int_I w(x_1,\dots,x_j,\dots,x_d) dx_j\)\(\frac1{\mes I} \int_I w^{-\frac1{p-1}}(x_1,\dots,x_j,\dots,x_d) dx_j\)^{{p-1}}\le C.
$$
In other words, belonging to ${\mathcal{A}}_p(\rd)$ implies
belonging to $A_p(\R^1)$ in each variable uniformly with respect to other variables;

\item[3)] For any $\delta=(\delta_1,\dots,\delta_d)\in \R^d_+$, $\delta_j>0$, $j=1,\dots,d$,
one has $w(\delta_1 x_1,\dots,\delta_d x_d)\in{\mathcal{A}}_p(\rd) $.

\end{enumerate}

\end{lem}

\begin{prop}
\label{prop1}
	Let  $1< q < \infty$, $\phi\in \cal B$, $w\in {\cal W}^{\alpha}$ for some $\alpha \in (0,1)$,
	and $w^q\in {\mathcal{A}}_q$. Then for any $f\in L_{q,w}$
	\begin{equation}\label{eqK1}
  \left(\sum_{k\in\zd}\left |{\langle f,{{\phi}_{0k}}\rangle}{w(k)}\right|^q\right)^{1/q}
	\le C \|f\|_{q,w},
\end{equation}
where $C$ depends only on $d$, $q$, $\alpha$, $c_w$ and $\phi$.
\end{prop}

 Before the proof of Proposition~\ref{prop1}, we introduce  additional notation and prove one lemma.
We set
$$
U_k^0=\{t\in \r\,:\,|t-k|<1\}\quad\text{and}\quad U_k^1=\r\setminus U_k^0,\quad k\in \z;
$$
if $k\in\zd$ and $\chi=(\chi_1,\dots,\chi_d)\in \{0,1\}^d$, then $U_k^\chi$ is defined by
$$
U_k^\chi = U_{k_1}^{\chi_1}\times\dots\times U_{k_d}^{\chi_d}.
$$
\begin{lem}\label{lemK1}
  Let $1<q<\infty$,  $w\in {\cal W}^{\alpha}(\R)$ for some $\alpha\in(0,1)$, and  $w^q\in A_q(\R)$. Then for any $f\in L_{q,w}(\r)$ and $u\in \r$
  \ban
	   \lll\sum_{k\in \z}\Bigg|{w(k)}\,\,\int\limits_{U_k^1}f(t)\frac{e^{2\pi i u(t-k)}}{t-k}dt\Bigg|^q\rrr^{1/q}
	\le C \|f\|_{L_{q,w}(\r)},
  \ean
	where $C$ depends only on $q$, $\alpha$, and $c_w$.
\end{lem}

{\bf Proof.}
Using properties of $w$, we have
\begin{equation}\label{26}
  \begin{split}
     \sum_{k\in \z}\Bigg|{w(k)}\,\,\int\limits_{U_k^1}f(t)\frac{e^{2\pi i u(t-k)}}{t-k}dt\Bigg|^q
&\le\sum_{k\in \z}\,\,\int\limits_{k-1/2}^{k+1/2}\,dx\Bigg|{w(x)}{w^*(x-k)}\,\,
\int\limits_{U_k^1}f(t)\frac{e^{2\pi i ut}}{t-k}dt\Bigg|^q\\
&\le  2^{\alpha q/2}c^q_w
\sum_{k\in \z}\,\,\int\limits_{k-1/2}^{k+1/2}\,dx\Bigg|{w(x)}\,\,\int\limits_{U_k^1}f(t)\frac{e^{2\pi i ut}}{t-k}dt\Bigg|^q.
  \end{split}
\end{equation}
Taking into account that by Minkowski's inequality
\begin{equation*}
  \begin{split}
      &\sum\limits_{k\,\in\,\z}\,\,\int\limits_{k-1/2}^{k+1/2}dx
    \( {w(x)}\int\limits_{U_k^1} {|f(t)|}
{\left|\frac{1}{t-x}-\frac{1}{t-k}\right|}\,dt\)^q\le
C_1
\sum\limits_{k\,\in\,\z}\,\,
\int\limits_{k-1/2}^{k+1/2}dx \( {w(x)}\!\!\int\limits_{|t-x|\ge 1/2}  \frac{|f(t)|}{(t-x)^{2}}\,dt\)^q
\\
 &\le C_1 \int\limits_{\r}dx
\( \int\limits_{\,\,|t|\ge 1/2}    \frac{c_w w^*(t)|w(t+x)f(t+x)|}
{t^{2}}\,dt\)^q\le
C_2 \|f\|^q_{q,w}\( \int\limits_{\,\,|t|\ge 1/2}    \frac{|t|^\alpha}
{t^{2}} \,dt\)^q=C_3\|f\|^q_{q,w},
   \end{split}
\end{equation*}
one can replace $t-k$ in the denominators of the integrand in~(\ref{26}) by $t-x$.

Next, we have
\begin{equation*}
  \begin{split}
     \int\limits_{\r} \bigg|
\int\limits_{|t-x|\ge1}f(t)
    \frac{e^{2\pi iut}} {t-x}\,dt\bigg|^q  \,{w^q(x)}\,dx
    &\le
    \int\limits_{\r}  \sup_{\epsilon>0}
    \bigg|\int\limits_{|t-x|\ge\epsilon}f(t)
 \frac{e^{2\pi iut}} {t-x}\,dt\bigg|^q\,{w^q(x)}\,dx\\
 &=\|\mathcal{M}(f_u)\|_{q,w}^q,
   \end{split}
\end{equation*}
where $f_u(t)=f(t)e^{2\pi iut}$ and $\mathcal{M}(g)$ is the maximal function of the
Hilbert transform of a function~$g$.   It remains to note that the operator $\mathcal{M}$
is bounded in $L_{q, w}$ under our assumption $w^q\in A_q(\R)$
(see, e.g.,~\cite[Corollary 7.13]{D}).
$\Diamond$

\bigskip

{\bf Proof of Proposition~\ref{prop1}.} For convenience, we introduce the following notation.
If $t\in\rd$, then  $\tilde{t}:=(t_1,\dots,t_{d-1})\in \R^{d-1}$.
For a function $g$ of $d$ variables $t_1,\dots, t_{d-1}$, $s$, we set $$g_s(\tilde{t})=g_{\tilde{t}}(s):=g(t_1,\dots,t_{d-1},s).$$
 Let also  $\tilde{S}=[a_1, b_1]\times\dots\times[a_{d-1}, b_{d-1}]$.

We have
$$
    \lll\sum_{k\in\Z^d}\left|{\langle f,\vp_{0k} \rangle}{w(k)}\right|^q\rrr^{1/q}
			=\lll\sum_{k\in\Z^d}\bigg|{w(k)}\int\limits_{\R^d} f(t)\vp(t-k)dt \bigg|^q\rrr^{1/q}
					\le \sum_{\chi\in\{0,1\}^d}\lll I^\chi\rrr^{1/q},
$$
where
$$
I^\chi=\sum_{k\in\Z^d}\bigg|w(k)\int\limits_{U_k^\chi} f(t)\vp(t-k)dt \bigg|^q.
$$
Thus, to prove~\eqref{eqK1} it suffices to show that for every $\chi\in \{0,1\}^d$
\begin{equation}\label{eqK9}
 I^\chi \le C_0 \Vert \t\Vert_{W_\infty^d(S)}^q \Vert f\Vert_{q,w}^q,
\end{equation}
where $C_0$ depends on $d$, $q$, $\alpha$, $S$, and $c_w$.

We  prove~\eqref{eqK9} by  induction on $d$. We will verify the  inductive step $d-1\to d$ and the base for $d=1$
simultaneously  using the same arguments.

To prove the inductive step $d-1\to d$, we assume that for any weight $\tilde{w}$ such that $\tilde{w}^q\in {\cal W}^{\alpha q}\cap \mathcal{A}_q({\R}^{d-1})$,
$g\in L_{q,\tilde{w}}(\R^{d-1})$, and $\tilde{\phi}\in \mathcal{B}(\R^{d-1})$
(more precisely $\tilde{\phi}=\mathcal{F}^{-1}\tilde{\t}$, where $\tilde{\t}$ is the same as in~\eqref{zzz61} with $\tilde{S}$ in place of $S$), we have
\begin{equation}\label{eqK10}
  \begin{split}
    \sum_{\tilde{k}\in\Z^{d-1}}  \bigg| {\tilde{w}(\tilde k)} \int\limits_{U_{\kt}^{\chit}} g\(\tt\)\tilde{\vp}\(\tt-\kt\)d\tt  \bigg|^q
		\le C_1\Vert \tilde{\t}\Vert_{W_\infty^{d-1}(\tilde{S})}^q\Vert g\Vert_{L_{q, \tilde{w}}(\R^{d-1})}^q,
  \end{split}
\end{equation}
where  $C_1$ depends on $d$, $p$, $\alpha$, $\tilde S$, and $C_{\tilde{w}}$.

In what follows,  we will use the fact  that under our assumptions, we have
$w^q_s\in{\cal W}^{\alpha q}\cap \mathcal{A}_q({\R}^{d-1})$
for every $s\in\r$ and
$w^q_{\tilde{t}}\in{\cal W}^{\alpha q}\cap A_q(\r)$
for every $\tilde{t}\in\r^{d-1}$. This follows from Lemma~\ref{lemAp} and basic properties of the weights in ${\cal W}^{\alpha q}$.

For any $\chi\in \{0,1\}^d$, we can write $\chi=(\chit,\chi_d)$ and
$$
I^\chi=\sum_{\chi_d\in \{0,1\}}I^{(\chit,\chi_d)}.
$$
Let us estimate $I^{(\chit,\chi_d)}$ for $\chi_d=0$ and $\chi_d=1$.

\smallskip

1)  First let $\chi_d=1$. In the case $d>1$, we set
\be
\label{300}
\psi(\tilde x, \eta) =\int\limits_{\tilde{S}}\theta_\eta(\tilde \xi) e^{2\pi i (\tilde\xi, \tilde{x})}\, d\tilde\xi.
\ee
Using the above notation and integrating by parts, we have
\begin{equation}\label{eqK7.5}
  \begin{split}
    \vp(x)&=\mathcal{F}^{-1}\t(x)=\int\limits_{a_d}^{b_d} \psi_{\tilde{x}}(\eta)e^{2\pi i x_d\eta}d\eta
		\\
    &=\frac{\psi_{\tilde{x}}(b_d)e^{2\pi ib_dx_d}-\psi_{\tilde{x}}(a_d)e^{2\pi ia_dx_d}}
		{2\pi ix_d}-\frac1{2\pi ix_d}\mathcal{F}^{-1}\psi'_{\tilde{x}}(x_d).
  \end{split}
\end{equation}
It follows from~\eqref{eqK7.5} that
\begin{equation}\label{eqK18}
  I^{(\chit,1)}=\sum_{\kt\in\Z^{d-1}}\sum_{l\in\Z}\bigg|
	\int\limits_{U_l^1}{ w_l(\tilde k)}\int\limits_{U_{\kt}^{\chit}} f_s\(\tt\){\vp}\(\tt-\kt,s-l\)d\tt ds\bigg|^q
\le (2\pi)^{-q}(I_1+I_2),
\end{equation}
where
\begin{equation*}
  I_1=\sum_{\kt\in\Z^{d-1}}\sum_{l\in\Z}\bigg| \int\limits_{U_l^1}
	{ w_l(\tilde k)}\int\limits_{U_{\kt}^{\chit}} f_s\(\tt\)\frac{\psi_{\tt-\kt}(b_d)e^{2\pi ib_d(s-l)}-
	\psi_{\tt-\kt}(a_d)e^{2\pi ia_d(s-l)}}{s-l}d\tt ds \bigg|^q
\end{equation*}
and
\begin{equation}\label{eqK19}
  I_2=\sum_{\kt\in\Z^{d-1}}\sum_{l\in\Z}\bigg| \int\limits_{U_l^1}{ w_l(\tilde k)}\int\limits_{U_{\kt}^{\chit}}
	f_s\(\tt\) \frac{\mathcal{F}^{-1}\psi'_{\tt-\kt}(s-l)}{s-l}d\tt ds \bigg|^q.
\end{equation}

Let us consider $I_1$.  Denoting
$$
  F_{\kt,u}(s)=\int\limits_{U_{\kt}^{\chit}} f_s(\tt)\psi_{\tt-\kt}(u) d\tt
$$
and using Lemma~\ref{lemK1}, we obtain
\begin{equation}\label{eqK19.5}
  \begin{split}
     I_1&\le \sum_{u\in\{a_d,b_d\}}\sum_{\kt\in\Z^{d-1}}
		\sum_{l\in\Z} \bigg|{w_{\tilde{k}}(l)}\int\limits_{U_l^1}\frac{F_{\kt,u}(s)e^{2\pi iu(s-l)}}{s-l}ds \bigg|^q
		\\
		&\le  C_2 \sum_{u\in\{a_d,b_d\}}\sum_{\kt\in\Z^{d-1}} \Vert F_{\kt,u}\Vert_{L_{q,  w_{\tilde{k}}}(\R)}^q,
  \end{split}
\end{equation}
where $C_2$ is the same constant  as in Lemma~\ref{lemK1}.

Now, using the induction hypothesis~\eqref{eqK10}, we derive
\begin{equation}\label{eqK19.6}
  \begin{split}
     \sum_{\kt\in\Z^{d-1}} \Vert F_{\kt,u}\Vert_{L_{q, w_{ \tilde{k}}}(\R)}^q
     &=\int\limits_\R\sum_{\kt\in\Z^{d-1}} \bigg|{ w_s(\tilde k)}
		\int\limits_{U_{\kt}^{\chit}} f_s(\tt)\mathcal{F}^{-1}\t_u(\tt-\kt) d\tt\, \bigg|^q \,{ds}
				\\
     &\le C_1 \Vert \t_u\Vert_{W_\infty^{d-1}(\widetilde{S})}^q
		\int\limits_\R \Vert f_s\Vert_{L_{q, w_s}(\R^{d-1})}^q \, {ds}.
  \end{split}
\end{equation}
Hence, combining~\eqref{eqK19.5} and~\eqref{eqK19.6}, we get
\begin{equation}\label{eqK20}
  I_1\le C_1 C_2 \sum_{u\in\{a_d,b_d\}}
	\Vert \t_u\Vert_{W_\infty^{d-1}(\widetilde{S})}^q \Vert f\Vert_{q,w}^q
	\le 2C_1C_2 \Vert \t\Vert_{W_\infty^{d}(S)}^q\Vert f\Vert_{q,w}^q.
\end{equation}
Let us consider $I_2$. Setting
$$
F_{\kt,\eta}^*(s)=\int\limits_{U_{\kt}^{\chit}} f_s(\tt)\psi_{\tt-\kt}'(\eta) d\tt
$$
and using H\"older's inequality, we obtain
\begin{equation}\label{eqK21}
  \begin{split}
    \Bigg|\int\limits_{U_l^1}\int\limits_{U_{\kt}^{\chit}} f_s(\tt)\frac{\mathcal{F}^{-1}\psi'_{\tt-\kt}(s-l)}{s-l}d\tt ds  \Bigg|^q
    &\le \(\int\limits_{a_d}^{b_d}
    \Bigg|\int\limits_{U_l^1}\int\limits_{U_{\kt}^{\chit}} f_s(\tt) \frac{\psi'_{\tt-\kt}(\eta)e^{2\pi i(s-l)\eta}}{s-l}d\tt ds\Bigg| d\eta\)^q
    \\
    &\le (b_d-a_d)^{q-1} \int\limits_{a_d}^{b_d}\Bigg|\int\limits_{U_l^1}F_{\kt,\eta}^*(s)\frac{e^{2\pi i(s-l)\eta}}{s-l} ds \Bigg|^qd\eta.
  \end{split}
\end{equation}
Thus, combining~\eqref{eqK19} and~\eqref{eqK21}, using Lemma~\ref{lemK1}, and the induction hypothesis~\eqref{eqK10}, we derive
\begin{equation}\label{eqK22}
  \begin{split}
    I_2 &\le (b_d-a_d)^{q-1}\sum_{\kt\in\Z^{d-1}}\sum_{l\in\Z} { w_l^q(\tilde k)} \int\limits_{a_d}^{b_d}
		\bigg|\int\limits_{U_l^1}F_{\kt,\eta}^*(s)\frac{e^{2\pi i(s-l)\eta}}{s-l} ds \bigg|^qd\eta
		\\
    &=(b_d-a_d)^{q-1}\int\limits_{a_d}^{b_d}\sum_{\kt\in\Z^{d-1}}\sum_{l\in\Z}
		\bigg|{ w_{\tilde k}(l)}\int\limits_{U_l^1}F_{\kt,\eta}^*(s)\frac{e^{2\pi i(s-l)\eta}}{s-l} ds \bigg|^qd\eta
		\\
    &\le  C_2(b_d-a_d)^{q-1}\int\limits_{a_d}^{b_d}\sum_{\kt\in\Z^{d-1}}
		\Vert F_{\kt,\eta}^* \Vert_{L_{q, w_{\tilde{k}}}(\R)}^q d\eta
		\\
    &=C_2(b_d-a_d)^{q-1}  \int\limits_{a_d}^{b_d}
		\int\limits_\R \sum_{\kt\in\Z^{d-1}}
		\bigg|{ w_s{(\tilde{k})}}\int\limits_{U_{\kt}^{\chit}} f_s(\tt)\mathcal{F}^{-1}
		\frac{\partial}{\partial \eta}\t_\eta\(\tt-\kt\)d\tt \bigg|^q \,{ds}\,d\eta
		\\
    &\le C_1C_2(b_d-a_d)^{q-1} \int\limits_{a_d}^{b_d}\bigg\Vert
		\frac{\partial}{\partial \eta}\t_\eta\bigg\Vert_{W_\infty^{d-1}(\widetilde{S})}^qd\eta
    \int\limits_\R \Vert f_s \Vert_{L_{q,w_s}(\R^{d-1})}^q\,{ds}
		\\
		&\le C_1C_2(b_d-a_d)^{q}\Vert \t\Vert_{W_\infty^{d}(S)}^q \Vert f\Vert_{L_{q,w}(\R^{d})}^q.
  \end{split}
\end{equation}

The above arguments are valid also for  $d=1$. In this case, the function $\psi_{\widetilde{x}}$ should be replaced by $\theta$ while $F_{\tilde k, u}(s)$ should be replaced by $f(s)\theta(u)$ in~\eqref{eqK19.5}.
Similarly,  $F^*_{\tilde k, \eta}(s)$ should be replaced by $f(s)\theta'(\eta)$ in~\eqref{eqK21}.
The sum over $\tilde k$ and the integral over ${U_{\kt}^{\chit}}$ are absent in this case.

Thus, combining~\eqref{eqK18}, \eqref{eqK20}, and~\eqref{eqK22}, we get~\eqref{eqK9} for any $d\ge 1$ and $\chi=(\chi_1,\dots,\chi_{d-1},1)\in \{0,1\}^d$.

\smallskip

2) Let now $\chi_d=0$.
In the case $d>1$, we have
\begin{equation*}
\begin{split}
  I^{(\chit,0)}&=\sum_{l\in\Z}\sum_{\kt\in\Z^{d-1}}\bigg|
	\int\limits_{U_l^0}{\tilde{w}_l(\tilde k)}
	\int\limits_{U_{\kt}^{\chit}} f_s\(\tt\){\vp}\(\tt-\kt,s-l\)d\tt ds\bigg|^q
	\\
&=\sum_{l\in\Z}
\sum_{\kt\in\Z^{d-1}}	\bigg|\int\limits_{U_l^0}{\tilde{w}_l(\tilde k)}
\int\limits_{U_{\kt}^{\chit}} f_s\(\tt\)\int\limits_{a_d}^{b_d}\,
\psi_\eta(\tilde t-\tilde k)e^{2\pi i\eta(s-l)}d\eta \,d\tt\,ds\bigg|^q,
	\end{split}
\end{equation*}
where $\psi$ is defined by~\eqref{300}.

Using H\"older's inequality and induction hypothesis~\eqref{eqK10}, we obtain
\begin{equation*}
  \begin{split}
      I^{(\chit,0)}&\le (b_d-a_d)^{q-1}
\sum_{l\in\Z}\,\int\limits_{|s-l|\le1}\int\limits_{a_d}^{b_d}
\sum_{\kt\in\Z^{d-1}}	\Bigg|{w_l(\tilde k)}
\int\limits_{U_{\kt}^{\chit}} f_s\(\tt\) \psi_\eta(\tilde t-\tilde k)e^{2\pi i\eta(s-l)} \,d\tt\Bigg|^q \,d\eta \,ds
\\
&\le C_1(b_d-a_d)^{q-1}
\sum_{l\in\Z}\,\int\limits_{|s-l|\le1}
\int\limits_{a_d}^{b_d}
\Vert \t_\eta\Vert_{W_\infty^{d-1}(\widetilde{S})}^q
\Vert f_s\Vert^q_{L_{q,  w_l}(\R^{d-1})}\,d\eta\,ds.
   \end{split}
\end{equation*}
Since
$$
w_l(\tilde{t})=w(t_1,\dots,t_{d-1},l)\le w^*(0,\dots,0,l-s)w(t_1,\dots,t_{d-1},s)
\le(1+|l-s|^2)^{\alpha/2}w_{s}(\tilde{t}),
$$
it follows that
\begin{equation*}
  \begin{split}
I^{(\chit,0)}&\le C_12^{q\alpha/2} (b_d-a_d)^q
\Vert \t\Vert_{W_\infty^{d}({S})}^q
\sum_{l\in\Z}\,\int\limits_{|s-l|\le1}
\Vert f_s\Vert^q_{L_{q, w_s}(\R^{d-1})}\,ds
\\
&\le C_{3} \Vert \t\Vert_{W_\infty^{d}({S})}^q
\Vert f\Vert^q_{L_{q, w}(\R^{d})},
   \end{split}
\end{equation*}
which yields~\eqref{eqK9} for any $\chi=(\chi_1,\dots,\chi_{d-1},0)\in \{0,1\}^d$.

If $d=1$, then due to properties of $w$ and H\"older's inequality, we have
\begin{equation*}\label{eqK3}
  \begin{split}
     \sum_{k\in \z} \Bigg|{w(k)}\int\limits_{U_k^0} f(t) \phi(t-k) dt\Bigg|^q
		&\le \|\phi\|_{L_\infty(\r)}^q \sum_{k\in\z}\lll \, \int\limits_{|t-k|\le1} w^*(t-k)\left|{f(t)}{w(t)}\right|dt\rrr^q
		\\
     &\le c_w^q 2^{\alpha q/2+q-1} (b_1-a_1)^q  \Vert\t\Vert_{L_\infty(S)}^q \sum_{k\in\z}\,
		\int\limits_{|t-k|\le1}\left|{f(t)}{w(t)}\right|^q dt\\
     &\le c_w^q 2^{\alpha q/2+q} (b_1-a_1)^q  \Vert\t\Vert_{L_\infty(S)}^q \|f\|_{L_{q,w}(\r)}^q,
  \end{split}
\end{equation*}
which yields~\eqref{eqK9} for $d=1$ and $\chi=0$.

\smallskip

Thus, inequality~\eqref{eqK9} holds for every $\chi\in \{0,1\}^d$ and $d\ge 1$, which proves the proposition.~$\Diamond$

\medskip

\begin{prop}
\label{coro2}
Let $1< p < \infty$, $\phi\in \cal B$,  $w\in {\cal W}^{\alpha}$ for some $\alpha \in (0,1)$, and $w^{-p}\in {\mathcal{A}}_p$. Then for any $a=\{a_k\}_{k\in \Z^d}\in {\ell_{p,1/w}}$
  \begin{equation}\label{221}
    \bigg\|\sum_{k\in\zd} a_k \phi_{0k}\bigg\|_{p,1/w}
		\le C \|a\|_{\ell_{p,1/w}},
  \end{equation}
	where $C$ depends on $d$, $p$, $\alpha$, $c_w$, and $\phi$.
 \end{prop}

{\bf Proof.}
 It follows from Riesz's theorem   that there exists $g\in L_{q,w}$,
  $1/p+1/q=1$, such that $\|g\|_{q, w}\le1$ and
$$
 \bigg\|\sum_{k\in\zd} a_k \phi_{0k}\bigg\|_{p,1/w}=
 \bigg|\sum_{k\in\zd} a_k\langle g,\phi_{0k}\rangle \bigg|=
\bigg|\sum_{k\in\zd}\frac{ a_k}{w(k)}\langle g,\phi_{0k}\rangle w(k) \bigg|.
$$
It is not difficult to
check that $w^{q}\in {\mathcal{A}}_q$.
Thus, to prove~\eqref{221} it remains to apply  H\"older's inequality and use Proposition~\ref{prop1} with $f=g$.\ $\Diamond$

\begin{prop}
\label{prop2}
	Let $1< p < \infty$, $\phi\in \cal B$,  $w\in {\cal W}^{\alpha}$ for some $\alpha \in (0,1)$, and $w^{-p}\in {\mathcal{A}}_p$. Then for any $f\in L_{p,1/w}$
	\begin{equation}\label{eqK1+}
  \left(\sum_{k\in\zd}\left |\frac{\langle f,{{\phi}_{0k}}\rangle}{w(k)}\right|^p\right)^\frac 1p
	\le C \|f\|_{p,1/w},
\end{equation}
where $C$ depends on $d$, $p$, $\alpha$, $c_w$ and $\phi$.
\end{prop}

{\bf Proof.} The proposition can be proved by following step by step the proof of Proposition~\ref{prop1} and using inequality~\eqref{1010} instead of~\eqref{W}.~$\Diamond$

\begin{coro}\label{coro3}
  Let $1< p<\infty$, $w\in {\cal W}^{\alpha}$ for some $\alpha \in (0,1)$, and $w^{-p}\in {\mathcal{A}}_p$. If  $\phi, \w \phi\in \cal B$, then for  any $f\in L_{p,1/w}$
  $$
    \bigg\|\sum_{k\in\zd} \langle f,\w\phi_{0k}\rangle \phi_{0k}\bigg\|_{p,1/w}
		\le C \Vert f\Vert_{p,1/w},
$$
	where $C$ depends on $d$,  $p$, $\alpha$, $c_w$, $\phi$, and $\w\phi$.
 \end{coro}

{\bf Proof.}
The proof follows immediately from Propositions~\ref{coro2} and~\ref{prop2}.~$\Diamond$

\begin{coro}\label{coro4}
  Let $1< p<\infty$, $w\in {\cal W}^{\alpha}$ for some $\alpha \in (0,1)$, and $w^{-p}\in {\mathcal{A}}_p$.
	If $\phi\in \cal B$ and  $\widetilde{\vp}\in \mathcal{L}_{q,w^*}$, $1/p+1/q=1$, then for  any $f\in L_{p,1/w}$
$$
    \bigg\|\sum_{k\in\zd} \langle f,\w\phi_{0k}\rangle \phi_{0k}\bigg\|_{p,1/w}
		\le C \Vert\w\vp\Vert_{\mathcal{L}_{q,w^*}} \Vert f\Vert_{p,1/w},
$$
		where $C$ depends on $d$, $p$, $\alpha$, $c_w$, and $\phi$.
 \end{coro}

{\bf Proof.}
The proof follows immediately from Propositions~\ref{prop1Lq} and ~\ref{coro2}.~$\Diamond$

\bigskip
\begin{lem}\label{lem1Band}
Let $1< p < \infty$,  $w\in {\cal W}^{\alpha}$ for some $\alpha \in (0,1)$,  $w^{-p}\in {\mathcal{A}}_p$, and $n\in \N$.
Suppose

\begin{enumerate}
  \item[1)]  $\psi\in\mathcal{B}$;
  \item[2)]  $\h{\w\psi}\in C^{\gamma}(\rd)$ for some $\gamma>n+d+p\alpha$ and $\h{\w\psi}$ is compactly supported;
  \item[3)]  $D^{\beta}\h{\w\psi}(\nul) = 0$ for all $[\beta]<n$, $\beta\in \zd_+$.
\end{enumerate}
Then for any $f\in W_{p,1/w}^n$
\begin{equation*}
   \bigg\|\sum\limits_{k\in\zd}
\langle f,\widetilde\psi_{0k}\rangle \psi_{0k}\bigg\|_{p, 1/w}\le  C \|f\|_{\dot W^n_{{p,1/w}}},
\end{equation*}
where $C$ depends only on $d$, $p$, $n$, $\psi$,  $\w \psi$, $c_w$, and $\alpha$.
\end{lem}

{\bf Proof.}
The lemma can be proved repeating step by step the proof of Lemma~\ref{lem1}. One needs only to use Proposition~\ref{coro2} instead of Proposition~\ref{uns2'}.~$\Diamond$

\bigskip

{\bf Proof of Theorem~\ref{lemJack+++}.} The proof is similar to the proof of Theorem~\ref{lemJack}. One only needs to use Corollary~\ref{coro4} instead of Corollary~\ref{lemBound}
\ $\Diamond$

\bigskip

{\bf Proof of Theorem~\ref{th4}.}
First we assume that $f\in W_{p,1/w}^n\cap L_2$ and prove that
\begin{equation}\label{101+}
   \bigg\|f-\sum\limits_{k\in\zd}
\langle f,\widetilde\phi_{0k}\rangle \phi_{0k}\bigg\|_{p, 1/w}\le  \Upsilon(w^*) \|f\|_{\dot W^n_{{p,1/w}}},
\end{equation}
where the functional $\Upsilon$ is independent of $f$ and
$
 \Upsilon(w^*(M^{-j}\cdot))\le C\Upsilon(w^*)
$
for all $j\in\z_+$, where the constant $C$ depends only on
$d$, $p$, $\phi$, $\w \phi$, $\alpha$, $C'$, and $c_w$.
To prove~\eqref{101+} we repeat step by step all arguments
of the proof of Lemma~\ref{lem2}, using Lemma~\ref{lem1Band}
instead of Lemma~\ref{lem1} in~\eqref{108} and
using Corollary~\ref{coro4}
instead of Corollary~\ref{lemBound} in~\eqref{108+}.
Next  we repeat  all steps of the proof of
 Theorem~\ref{th2}, using~\eqref{101+} instead of Lemma~\ref{lem2} and
 taking into account that  $(w(M^{-j}\cdot))^{-p}\in  {\mathcal{A}}_p$ by Lemma~\ref{lemAp} and $\Vert M^{-j}\Vert=|\lambda|^{-j}$.
\ $\Diamond$

\subsection{The case of sampling expansions}\label{S5.3}
In this section, we study expansions
$\sum\limits_{k\in\zd}
\langle f,\widetilde\phi_{jk}\rangle \phi_{jk}$
 where  $\w\vp$ is the Dirac delta-function, i.e.,
$$
\sum\limits_{k\in\zd}
\langle f,\widetilde\phi_{jk}\rangle \phi_{jk}=
\sum\limits_{k\in\zd}
\langle \h f, \widehat{\delta_{jk}}\rangle \phi_{jk}=
m^{-j/2}\sum\limits_{k\in\zd} f(- M^{-j}k) \phi_{jk}.
$$

In the theorems below, we will suppose that the function $\vp$ and the Dirac delta-function are strictly compatible (weakly compatible of order $n\in \N$), which implies that
$\widehat \phi(\xi)\equiv1$ on the ball $B_{\delta}$  ($D^{\beta}(1-\h\phi)({\bf 0}) = 0$ for all $[\beta]<n$, $\beta\in\zd_+$).

First, we consider band-limited weights $w$ and, as in the case of Theorem~\ref{th4}, we suppose that  $w\in {\cal W}^{\alpha}$ and $w^{-p}\in {\mathcal{A}}_p$.

\begin{theo}
\label{theoQj_new}
	Let  $2\le p < \infty$, $w\in {\cal W}^{\alpha}$ for some $\alpha \in (0,1)$,
	$w^{-p}\in {\mathcal{A}}_p$, $\supp \widehat{w}\subset B_{\delta/2}$ for some  $\delta\in (0,1/2)$, and let  $M \in \mathfrak{M}$ be a diagonal  matrix.
	Suppose

\begin{enumerate}
  \item[1)]  $\phi\in \cal B$;

   \item[2)]  $\phi$ is strictly compatible with the Dirac delta-function with respect to the parameter $\delta$;

  \item[3)]  $\supp \h\vp \subset (-1+\delta/2,1-\delta/2)^d$.
\end{enumerate}
If a function $f\in L_{p,1/w}$ is such that
$\mathcal{F}({w^{-1}f})\in L_q$, $1/p+1/q=1$, and $\mathcal{F}({w^{-1}f})(\xi)=
\mathcal{O}(|\xi|^{-d-a})$ as $|\xi|\to\infty$, $a>0$, then for any $j\in \Z_+$ and $\gamma> d/p$
\be
\bigg\Vert f- m^{-j/2}\sum\limits_{k\in\zd} f(- M^{-j}k) \phi_{jk}\bigg\Vert_{p, 1/w}^q
\le
C \|M^{-j}\|^{q\gamma}\!\!  \int\limits_{|M^{-j}\xi|\ge\delta/2}
|\xi|^{q\gamma}| \mathcal{F}({w^{-1}f})(\xi)|^q d\xi,
\label{8+}
\ee
where $C$ does not depend on  $j$ and~$f$.
\end{theo}	

{\bf Proof of Theorem~\ref{theoQj_new}.}
First, we consider the case $j=0$ and prove that
$$
	\bigg\Vert f- \sum\limits_{k\in\zd}
 f(- k) \phi_{0k}\bigg\Vert_{p, 1/w}^q\le C_1
	 \int\limits_{|\xi|\ge\delta/2}
	|\xi|^{q\gamma}| \mathcal{F}({w^{-1}f})(\xi)|^q d\xi,
$$
where $C_1$ depends on  $p$, $\alpha$, $c_w$ and $\phi$.

We set
$$
	 g=w^{-1}f\quad\text{and}\quad G(\xi)=\sum\limits_{l\in\,\zd}
	\h g(\xi+l).
$$
 Using  Lemma 1 from \cite{KS1}  and the  Hausdorff-Young inequality, we have
$$
\langle \h g,\h{{\delta}_{0k}}\rangle=\h G(k)
$$
and
	\be
	\left(\sum_{k\in\zd} |\langle \h g,\h{{\delta}_{0k}}\rangle|^p\right)^{1/p} =
	\left(\sum_{k\in\zd} |\h G(k)|^p\right)^{1/p} \le
	\|G\|_{L_q(\td)} <\infty.
	\label{1_prop}
	\ee
Since
$$
\langle \h g,\h{{\delta}_{0k}}\rangle=\langle \mathcal{F}({w^{-1}f}),\h{{\delta}_{0k}}\rangle=(w^{-1}f)(-k),
$$
it follows from Proposition~\ref{coro2} and inequality~\eqref{1_prop} that
	\be
		\bigg\| \sum\limits_{k\in\zd}
 f(- k) \phi_{0k}\bigg\|_{p, 1/w}\le C_2
	\left(\sum_{k\in\zd}
	\bigg|\frac{f(-k)}{w(-k)}\bigg|^p\right)^{1/p}
	\le C_3 \|G\|_{L_q(\td)}.
	\label{0_prop}
	\ee
Again using  Lemma 1 from \cite{KS1}, we have
	\be
	\label{209}
\|G\|_{L_q(\td)}
	\le 	C_4\lll\( \,\int\limits_{|\xi|\ge\delta}
	|\xi|^{q\gamma}| \h g(\xi)|^q d\xi\)^{1/q}+\|\h g\|_q\rrr.
	\ee
Due to the du Bois-Reymond lemma,  the function $g(-x)$ coincides with  $\h {\h g}(x)$ almost everywhere.
It follows from  the Hausdorff-Young inequality that
$\|g\|_p\le\|\h g\|_q$, which together with~\eqref{0_prop} and \eqref{209} yields
\begin{equation}\label{35}
  \begin{split}
&\bigg\|f- \sum\limits_{k\in\zd} f(- k) \phi_{0k}\bigg\|_{p, 1/w}
\le\|g\|_p+\bigg\| \sum\limits_{k\in\zd} f(- k) \phi_{0k}\bigg\|_{p, 1/w}
\\
&\le	C_5\lll \,\int\limits_{|\xi|\ge\delta} |\xi|^{q\gamma}| \h g(\xi)|^q d\xi+\int\limits_{|\xi|<\delta}| \h g(\xi)|^qd\xi\rrr^{1/q}
\\
&=C_5 \lll \,\int\limits_{|\xi|\ge \delta}|\xi|^{\gamma q}| \mathcal{F}({w^{-1}f})(\xi)|^q d\xi+\int\limits_{|\xi|<\delta}| \mathcal{F}({w^{-1}f})(\xi)|^qd\xi\rrr^{1/q}.
  \end{split}
\end{equation}
		
Denote
$$
F(x)=w(x)H(x),\quad H(x)=\mathcal{F}^{-1} h(x),
$$
where $h\in  C^\infty(\R^d)$ and $\supp h\subset B_{\delta/2}$. We can choose the function $h$ such that
\begin{equation}\label{chet}
  \int_{|\xi|<\delta/2}|\mathcal{F}({w^{-1}f})(\xi)-h(\xi)|^q d\xi<\int_{|\xi|>\delta}|\xi|^{\gamma q}|\mathcal{F}({w^{-1}f})(\xi)|^q d\xi.
\end{equation}

Note that $F\in L_2$ since $|F(x)|\le w(x)|H(x)|$ and $|H(x)|\le c_N (1+|x|)^{-N}$ for any $N\in \N$. Moreover, by the Paley-Wiener-Schwartz theorem, $\supp \h F\subset B_\delta$. Thus, by Proposition~\ref{corE}, \eqref{35}, and~\eqref{chet}, we obtain
\begin{equation*}
  \begin{split}
        &\bigg\|f- \sum\limits_{k\in\zd}
 f(- k) \phi_{0k}\bigg\|_{p, 1/w}^q=\bigg\|f-F- \sum\limits_{k\in\zd}
 (f(- k)-F(-k)) \phi_{0k}\bigg\|_{p, 1/w}^q\\
 &\le C_5^q\lll \,\,\int\limits_{|\xi|\ge \delta}|\xi|^{\gamma q}| \mathcal{F}\({w^{-1}(f-F)}\)(\xi)|^q d\xi+\int\limits_{|\xi|<\delta}| \mathcal{F}\({w^{-1}(f-F)}\)(\xi)|^qd\xi\rrr\\
&=C_6^q\int\limits_{|\xi|\ge \delta/2}|\xi|^{\gamma q}|
\mathcal{F}\({w^{-1}f}\)(\xi)|^q d\xi.
  \end{split}
\end{equation*}

This yields~(\ref{8+}) for $j=0$. To get the required inequality for any $j\in \N$,
we replace $f$ by $f(M^{-j} \cdot)$ and $w$ by $w(M^{-j} \cdot)$,
 take into account that  $(w(M^{-j}\cdot))^{-p}\in  {\mathcal{A}}_p$
by Lemma~\ref{lemAp}, and change variables in the integrals. $\Diamond$	

\begin{coro}\label{cordecanalweigh}
Under the assumptions of Theorem~\ref{theoQj_new}, we have
$$
   \bigg\Vert f- m^{-j/2}\sum\limits_{k\in\zd} f(- M^{-j}k) \phi_{jk}\bigg\Vert_{p, 1/w}= \mathcal{O}\(|\lambda|^{-j(a+d/p)}\),
$$
where $\lambda$ is the smallest (in absolute value) diagonal element of $M$.
\end{coro}

\begin{rem}
In Theorem~\ref{theoQj_new} and Corollary~\ref{cordecanalweigh}, we suppose that $\supp \widehat{w}\subset B_{\delta/2}$,  which is quite exotic in such type of problems.
Nevertheless, such weights can be easily  constructed in the following way.
For some  weight $v$ such that $v\in {\cal W}^{\alpha}$ and $v^{-p}\in {\mathcal{A}}_p$, we set
$$
w(x)=v*V(x),
$$
where $V$ is such that $\supp \h V$ is compact and $V(x)\ge 0$ for all $x\in \R^d$.
It is not difficult to see that if $v$ is positive, symmetric, and $v(x+y)\le  v^*(x)v(y)$, then  for all $x\in \R^d$ one has the following two-sided inequality:
\begin{equation*}\label{twosided}
  v(x)\Vert V \Vert_{L_{1,1/v^*}}\le w(x)\le v(x)\Vert V \Vert_{L_{1,v^*}}.
\end{equation*}
Thus, for an appropriate function $V$, we have that  $w\in {\cal W}^{\alpha}$ with $w^*=v^*$, $w^{-p}\in {\mathcal{A}}_p$, and $\supp \h w$ is compact.
\end{rem}


\begin{theo}
\label{theoQj_new1}
	Let  $2\le p < \infty$, $n\in \N$, $w\in {\cal W}^{\alpha}$
	for some $\alpha \in (0,1)$, $w^{-p}\in {\mathcal{A}}_p$,
	$w\in C^\infty(\rd)$, $|D^\beta w(x)|\le c_{w,\beta} w(x)$
	for every $\beta\in\zd_+$, $x\in \R^d$, and let  $M \in \mathfrak{M}$ be a diagonal  matrix.
	Suppose
\begin{enumerate}
  \item[1)] $\phi\in \cal B$;
  \item[2)] $\h\phi \in C^{r}(B_{\varepsilon})$ for some integer $r>n+d+p\alpha$ and $\varepsilon>0$;
  \item[3)] $\phi$ is weakly compatible of order $n$  with the Dirac delta-function;
  \item[4)]  $\supp \h\vp \subset (-1,1)^d$.
\end{enumerate}
If a function $f\in L_{p,1/w}$ is such that
$\mathcal{F}({w^{-1}f})\in L_q$, $1/p+1/q=1$, and $\mathcal{F}({w^{-1}f})(\xi)=
\mathcal{O}(|\xi|^{-d-a})$ as $|\xi|\to\infty$, $a>0$, then for any $j\in \Z_+$ and any $\gamma> d/p$
\be
\begin{split}
   \bigg\Vert f- m^{-j/2}\sum\limits_{k\in\zd} f(- M^{-j}k) \phi_{jk}\bigg\Vert_{p, 1/w}^q
&\le C_1 \sum_{\nu=0}^n \|M^{-j}\|^{\nu q} \omega_{n-\nu}(f,\|M^{-j}\|)_{p,1/w}^q\\
&\qquad\qquad+C_2 \|M^{-j}\|^{q\gamma }\!\!  \int\limits_{|M^{-j}\xi|\ge 1/2}
|\xi|^{q\gamma}\Big| \mathcal{F}{\Big(\frac fw\Big)}(\xi)\Big|^q d\xi,
\end{split}
\label{8++++}
\ee
where $C_1$ and $C_2$ do not depend on $j$ and $f$.
\end{theo}

\begin{rem}
Theorem~\ref{theoQj_new1} is new also in the unweighed case, i.e. for $w(x)\equiv 1$.
In this case, one can show that inequality~\eqref{8++++} holds for any  $M \in \mathfrak{M}$
and has the following form:
\begin{equation*}
  \begin{split}
      \bigg\Vert f- m^{-j/2}\sum\limits_{k\in\zd} f(- M^{-j}k) \phi_{jk}\bigg\Vert_{p}^q
&\le C_1 \omega_{n}(f,\|M^{-j}\|)_{p}^q\\
&\qquad\qquad+C_2 \|M^{-j}\|^{q\gamma }\!\!  \int\limits_{|M^{-j}\xi|\ge 1/2}
|\xi|^{q\gamma}| \mathcal{F}{f}(\xi)|^q d\xi,
   \end{split}
\end{equation*}
where $C_1$ and $C_2$ do not depend on $j$ and $f$.
\end{rem}

\begin{rem}
\label{rem1++_}
Similarly to Theorem~\ref{th4} and Theorem~\ref{th2}, it is possible to obtain a sharper  version of inequality~\eqref{8++++} by
replacing the modulus of smoothness $\omega_{n-\nu}\(f, \|M^{-j}\|\)_{p,1/w}$ by the anisotropic modulus of smoothness $\Omega_{n-\nu}\(f, M^{-j}\)_{p,1/w}$.
\end{rem}

{\bf Proof of Theorem~\ref{theoQj_new1}.}
Let $\w\phi$ be a function such that  $\h{\w\phi}$ is infinitely differentiable, $\supp \widehat{\w\phi}\subset [-1, 1]^d$,  and
$ \h{\w\phi}(\xi)\equiv 1$ on $\T^d$. Since
$$
f(-k)=w(-k)\Big\langle \mathcal{F}{\Big(\frac fw\Big)}, \h{\delta_{0k}}\Big\rangle,
$$
we have
\ba
\label{502}
\begin{split}
    &\Big\|f-\sum\limits_{k\in\zd}f(-k)\phi_{0k}\Big\|_{p, 1/w}\\
    &\le
\Big\|\sum\limits_{k\in\zd}w(-k)\Big\langle \mathcal{F}{\Big(\frac fw\Big)},
\h{\delta_{0k}}-\h{\w\phi_{0k}}\Big\rangle\phi_{0k}\Big\|_{p, 1/w}
+ \Big\|f-\sum\limits_{k\in\zd}w(-k)\Big\langle {\frac fw},
{\w\phi_{0k}}\Big\rangle\phi_{0k}\Big\|_{p, 1/w}\\
&=:I_1(f,w)+I_2(f,w).
\end{split}
\ea
Using Proposition~\ref{coro2}, the Hausdorff-Young inequality, Lemma~1 in \cite{KS1}, and
taking into account that $ 1-\h{\w\phi}(\xi)=0$ if $\xi\in \T^d$, we obtain
\begin{equation}\label{503}
  \begin{split}
I_1(f,w) &\le C_1\lll
\sum\limits_{k\in\zd}\Big|\Big\langle \mathcal{F}{\Big(\frac fw\Big)},
\h{\delta_{0k}}-\h{\w\phi_{0k}}\Big\rangle\Big|^p\rrr^{1/p}\\
&\le
C_2\Bigg\|\sum_{l\in\zd}\mathcal{F}{\Big(\frac fw\Big)}(\cdot+l)\(1-\overline{\h{\w\phi}(\cdot+l)}\) \Bigg\|_{L_q\(\T^d\)}\\
&\le
C_3\lll\,\,\int\limits_{|\xi|\ge 1/2}|\xi|^{q\gamma}
\bigg|\mathcal{F}{\Big(\frac fw\Big)}(\xi)\bigg|^q\, d\xi\rrr^{1/q}.
   \end{split}
\end{equation}
The functions $\phi, \w\phi$ satisfy all assumptions of Theorem~\ref{th4}, hence,
\begin{equation}\label{504}
  \begin{split}
     I_2(f,w)&\le \bigg\|f-\sum\limits_{k\in\zd}\langle f,
{\w\phi_{0k}}\rangle\phi_{0k}-\sum\limits_{k\in\zd}\Big\langle {\frac {w(-k)-w}w}f,
{\w\phi_{0k}}\Big\rangle\phi_{0k}\bigg\|_{p, 1/w}
\\
&\le C_4 \omega_{n}\(f,1\)_{p,1/w}+ I_3(f,w),
  \end{split}
\end{equation}
where
$$
I_3(f,w)=\bigg\|\sum\limits_{k\in\zd}\int\limits_{\rd} {\frac {w(-k)-w(t)}{w(t)}f(t)
\overline{\w\phi_{0k}}(t)}\, dt\,\phi_{0k}\bigg\|_{p, 1/w}.
$$

Using Taylor's formula, we have
$$
w(-k)-w(t)=\sum_{0<[\beta]\le n}\frac{(-1)^{[\beta]}D^\beta w(t)}{\beta!} (t+k)^\beta+
\sum_{[\beta]=n+1}\frac{(-1)^{[\beta]}D^\beta w(t+\eta(t+k))}{\beta!} (t+k)^\beta
$$
for some $\eta\in(-1,0)$, which gives
\begin{equation*}
  I_3(f,w)\le \sum_{0<[\beta]\le n} I_{3,\beta}(f,w)+\sum_{[\beta]=n+1} I_{4,\beta}(f,w),
\end{equation*}
where
$$
I_{3,\beta}(f,w)=\bigg\|\sum\limits_{k\in\zd}\Big\langle{\frac {D^\beta w}{w}f,
{\psi_{0k}}}\Big\rangle\phi_{0k}\bigg\|_{p, 1/w},
$$
$$
I_{4,\beta}(f,w)=\bigg\|\sum\limits_{k\in\zd}\Big\langle{\frac {D^\beta w(\cdot+\eta(\cdot+k))}{w}f,
{\psi_{0k}}}\Big\rangle\phi_{0k}\bigg\|_{p, 1/w},
$$
and
$
\psi(t)=t^\beta\w\phi(t).
$

Fix $\beta\in\zd_+$, $0<[\beta]\le n$.  Let a function $g$ satisfy~\eqref{Kf1} and~\eqref{Kf2} with $\frac{D^\beta w}{w}f$ instead of $f$. Using Corollary~\ref{coro4} and Lemma~\ref{lem1Band}, taking into account that $\h\psi(\xi)=(-2\pi i)^{-[\beta]}D^\beta\h{\w\phi}(\xi)$,
and hence $D^\beta \h\psi(\nul)=0$ for any $\beta\in\zd_+\setminus\{\nul\}$, we obtain
\begin{equation}\label{newFFF}
  \begin{split}
      I_{3,\beta}(f,w)&\le \bigg\|\sum\limits_{k\in\zd}\Big\langle{\frac {D^\beta w}{w}f-g,
{\psi_{0k}}}\Big\rangle\phi_{0k}\bigg\|_{p, 1/w}+\bigg\|\sum\limits_{k\in\zd}\langle{g,
{\psi_{0k}}}\rangle\phi_{0k}\bigg\|_{p, 1/w}\\
&\le C_5\(\bigg\|\frac{D^\beta w}{w}f-g\bigg\|_{p, 1/w}+\Vert g\Vert_{\dot W_{p,1/w}^n} \)\le C_6\omega_n\(\frac{D^\beta w}{w}f,1\)_{p,1/w}.
  \end{split}
\end{equation}
Combining~(\ref{502}), (\ref{503}), (\ref{504}), and~\eqref{newFFF}, we derive
\begin{equation*}
  \begin{split}
     &\Big\|f-\sum\limits_{k\in\zd}f(-k)\phi_{0k}\Big\|_{p, 1/w}
\le C_7 \lll\,\,\int\limits_{|\xi|\ge 1/2}|\xi|^{q\gamma}
\bigg|\mathcal{F}{\Big(\frac fw\Big)}(\xi)\bigg|^q\, d\xi\rrr^{1/q}
\\
&+C_{8} \omega_{n}\(f,1\)_{p,1/w}+C_{9}\sum_{0<[\beta]\le n}\omega_n\(\frac{D^\beta w}{w}f,1\)_{p,1/w}+C_{10}\sum_{[\beta]=n+1} I_{4,\beta}(f,w).
  \end{split}
\end{equation*}
For an arbitrary $j\in\z_+$, taking into account that  $(w(M^{-j}\cdot))^{-p}\in  {\mathcal{A}}_p$,
by Lemma~\ref{lemAp}, we can replace $f$ by $f(M^{-j} \cdot)$ and  $w$ by $w(M^{-j} \cdot)$,
  and after  change of variables in all integrals, we have
\begin{equation}\label{DOPPP-}
  \begin{split}
&\bigg\|f- m^{-j/2}\sum\limits_{k\in\zd} f(- M^{-j}k) \phi_{jk}\bigg\|_{p, 1/w}\\
&=m^{-j/p} \bigg\|f(M^{-j}\cdot)-\sum\limits_{k\in\zd} f(M^{-j}k) \phi_{0k}\bigg\|_{p,1/w(M^{-j}\cdot)}
\\
&\le C_{7}\|M^{-j}\|^{\gamma }\lll
\,\int\limits_{|M^{-j}\xi|\ge 1/2}
|\xi|^{q\gamma}| \mathcal{F}({w^{-1}f})(\xi)|^q d\xi
\rrr^{1/q}+ C_{8} \omega_{n}\(f,\|M^{-j}\|\)_{p,1/w}
\\
&+C_{9}\sum_{0<[\beta]\le n}\omega_n\(\frac{D^\beta w}{w}f,\|M^{-j}\|\)_{p,1/w}
+C_{10}\sum_{[\beta]=n+1} m^{-j/p}I_{4,\beta}(f_j,w_j),
   \end{split}
\end{equation}
where $f_j(t)=f(M^{-j}t)$, $w_j(t)=w(M^{-j}t)$.

To estimate the modulus $\omega_n(\frac{D^\beta w}{w}f,\|M^{-j}\|)_{p,1/w}$, we use the following well-known relations:
\be
\label{1000}
\Delta_h^n (f_1 f_2)=\sum_{\nu=0}^n \binom{n}{\nu} \Delta_h^\nu (f_1) \Delta_h^{n-\nu}(f_2),
\ee
$$
\omega_\nu (f_1,h)_\infty \le C(\nu) h^\nu \Vert f_1\Vert_{\dot W_\infty^\nu},\quad f_1\in W_{\infty}^\nu,
$$
which imply that
\begin{equation}\label{DOPPP}
  \begin{split}
    \omega_n\(\frac{D^\beta w}{w}f,\|M^{-j}\|\)_{p,1/w}&
		\le \sum_{\nu=0}^n \binom{n}{\nu}\omega_\nu\(\frac{D^\beta w}{w},\|M^{-j}\|\)_{\infty}
    \omega_{n-\nu}\(f,\|M^{-j}\|\)_{p,1/w}\\
    &\le C_{11}\sum_{\nu=0}^n  \|M^{-j}\|^{\nu } \omega_{n-\nu}(f,\|M^{-j}\|)_{p,1/w}.
  \end{split}
\end{equation}
To estimate $I_{4,\beta}(f_j, w_j)$ for $[\beta]=n+1$, we note that
\begin{equation*}
  \begin{split}
     |D^\beta w_j(t+\eta(t+k))|&\le
\|M^{-j}\|^{n+1} \left| D^\beta w\(M^{-j}t+\eta M^{-j}(t+k)\)\right|\\
&\le c_{w,\beta} \|M^{-j}\|^{n+1} w(M^{-j}t+\eta M^{-j}(t+k))\\
&\le c_{w,\beta}c_w \|M^{-j}\|^{n+1}w_j(t) (1+|t+k|^\alpha).
   \end{split}
\end{equation*}

Thus, using Proposition~\ref{coro2} and Proposition~\ref{prop1Lq},
taking into account that the  function
 $(1+|t|^{\alpha})|t^{\beta}\w\phi(t)|$ belongs to ${\mathcal{L}_{q,w^*}}$, we obtain
  \begin{equation*}
  \begin{split}
       I_{4,\beta}(f_j, w_j)
       &\le c_{w,\beta}c_w \|M^{-j}\|^{n+1}
\left\|\sum\limits_{k\in\zd}\int\limits_{\rd} |f_j(t)|
(1+|t+k|^{\alpha})|(t+k)^\beta\w\phi_{0k}(t)|\, dt\,\phi_{0k}\right\|_{p, 1/w_j}
\\
&\le C_{12}\|M^{-j}\|^{n+1}\(\sum_{k\in\Z^d} \Bigg(\frac1{w_j(k)} \int\limits_{\rd} |f_j(t)|
(1+|t+k|^{\alpha})|(t+k)^\beta\w\phi(t+k)|dt \Bigg)^p\)^{1/p}
\\
&\le C_{13}\|M^{-j}\|^{n+1} \Vert f_j\Vert_{p,1/w_j}=C_{13}\|M^{-j}\|^{n+1} m^{j/p}  \Vert f\Vert_{p,1/w}.
   \end{split}
\end{equation*}
Combining this with~\eqref{DOPPP-} and \eqref{DOPPP},  we get~\eqref{8++++}.~$\Diamond$

\begin{rem}
\label{rem1+++_++}
 Analysing the proof of Theorem~\ref{theoQj_new1}, it is clear that inequality~\eqref{8++++} remains valid for any $2\le p\le \infty$ under the assumption  $\phi\in \mathcal{L}_{p,w^*}\cap L_2$ (instead of $\phi\in \mathcal{B}$). In this case, the condition $w^{-p}\in {\mathcal{A}}_p$ can be dropped and it suffices to assume that $\alpha>0$.
\end{rem}	

\begin{rem}
\label{rem1+++_}
The first term in the right hand side of~\eqref{8++++} can be replaced by
$$
C_1 \sum_{\nu=0}^n \|M^{-j}\|^{\nu q} \omega_{n-\nu}\(\frac f w,\|M^{-j}\|\)_{p}^q.
$$
Indeed, it follows from~\eqref{1000} that
 for any $\mu=1,\dots, n$ and $\delta=\|M^{-j}\|$
$$
\omega_{\mu}\( f, \delta\)_{p, 1/w}=\sup\limits_{|t|\le \delta}
\Big\|\frac{\Delta_t^\mu f}w\Big\|_p\le\sum_{\nu=0}^\mu \binom{\mu}{\nu} \sup\limits_{|t|\le \delta}\Big\|\frac{ \Delta_t^{\nu}w}w\Big\|_\infty\omega_{\mu-\nu} \(\frac f w, \delta\)_p.
$$
Using Taylor's formula with Lagrange's remainder, we have
$$
\Big|\frac{ \Delta_t^{\nu}w(x)}{w(x)}\Big|\le C(n)\sum_{l=1}^\nu\sum_{[\beta]=
 \nu}\Big| w^{-1}(x)D^\beta w(x+\theta_{l} t)\Big||t|^\nu,
$$
where $|\theta_{l}|\le \nu$. It remains to note that, due to properties of $w$,
$$
\Big|\frac{D^ \beta w(x+\theta_{l} t)}{w(x)}\Big|\le w^*(\theta_{l} t)
\Big|\frac{D^ \beta w(x+\theta_{l} t)}{w(x+\theta_{l} l)}\Big|\le C(w,n)<\infty.
$$
\end{rem}

\begin{coro}\label{coro22}
Under the assumptions of Theorem~\ref{theoQj_new1}, the following estimate holds for $j\to \infty$
$$
\bigg\Vert f- m^{-j/2}\sum\limits_{k\in\zd} f(- M^{-j}k) \phi_{jk}\bigg\Vert_{p, 1/w}=
 \begin{cases}
	 	\mathcal{O}(|\lambda|^{-j(d/p+a)})  &\mbox{if } n>d/p+a,
\\
	 	\mathcal{O}(|\lambda|^{-jn}j^{1/2})  &\mbox{if } n=d/p+a,
\\
\mathcal{O}(|\lambda|^{-jn})  &\mbox{if }  n<d/p+a,
	 	 \\
	\end{cases}
$$
where $\lambda$ is the smallest (in absolute value) diagonal element of $M$.
\end{coro}

{\bf Proof.}
Obviously,
$$
\|M^{-j}\|^{q\gamma }\!\!  \int\limits_{|M^{-j}\xi|\ge\delta}
|\xi|^{q\gamma}\Big| \mathcal{F}{\Big(\frac fw\Big)}(\xi)\Big|^q d\xi=
\mathcal{O}\(|\lambda|^{-jq(a+d/p)}\),\quad j\to \infty.
$$
It remains to estimate the first term in the right hand side of~\eqref{8++++}.
Due to Remark~\ref{rem1+++_}, it suffices to estimate the sum $\sum_{\nu=0}^n \delta^{\nu q} \omega_{n-\nu}\(w^{-1}f,\delta\)_{p}^q$.

Set
$$
V_\sigma f(x)=\mathcal{F}^{-1}\(v\(\sigma^{-1}|\xi|\)\widehat{f}(\xi)\)(x),
$$
where $v\in C^\infty(\R)$, $v(\xi)\le1$,  $v(\xi)=1$ for $|\xi|\le 1/2$ and $v(\xi)=0$ for $|\xi|\ge 1$. Using
Pitt's inequality (see, e.g., \cite[inequality (1.1) for $s=0$ and $p=q$]{GT}), we have
$$
\Vert \widehat{g}\Vert_{L_p(\R^d)}\le C(p)\(\,\int\limits_{\R^d} |x|^{d(p-2)}|g(x)|^p dx\)^{1/p},\quad 2<p<\infty.
$$
It follows  that
\begin{equation*}
  \begin{split}
      E_{B_\sigma}(w^{-1}f)_p&\le \Vert w^{-1}f-V_\sigma(w^{-1}f)\Vert_p\\
      &\le C(p)\(\,\int\limits_{|\xi|\ge \sigma/2}|\xi|^{d(p-2)} |(1-v(\sigma^{-1}|\xi|))\mathcal{F}(w^{-1}f)(\xi)|^p d\xi\)^{1/p}\\
      &=\mathcal{O}\(\sigma^{-(d/p+a)}\),\quad \sigma\to \infty.
   \end{split}
\end{equation*}
If $n\ge d/p+a$, using the following Marchaud inequality (see~\cite{DDT}):
\begin{equation}\label{eqMarch}
  \omega_\nu(w^{-1}f,2^{-N})_p\le C2^{-\nu N}\(\sum_{k=0}^N 2^{2\nu k} E_{B_{2^k}}(f)_p^2\)^{1/2},
\end{equation}
we obtain
$$
\(\sum_{\nu=0}^n \delta^{\nu q} \omega_{n-\nu}\(w^{-1}f,\delta\)_{p}^q\)^\frac1q = \left\{
                                                                         \begin{array}{ll}
                                                                           \mathcal{O}\(\delta^{a+d/p}\), & \hbox{$n> d/p+a$,} \\
                                                                           \mathcal{O}\(\delta^{n}\log^{1/2}(\delta^{-1})\), & \hbox{$n=d/p+a$.}
                                                                         \end{array}
                                                                       \right.
$$
If now $n<a+d/p$, then $d+a>n+d/q$, and hence the function
$|\xi|^{n} \mathcal{F}{(w^{-1} f)}(\xi) $ belongs to $L_q$.
 It follows from  the Hausdorff-Young inequality that $w^{-1}f\in W^n_p$, which yields that
$$
\omega_\nu(w^{-1}f,\delta)_p=\mathcal{O}(\delta^{\nu}),\quad \delta\to 0,
$$
for any $1\le \nu\le n$, which proves the corollary. \ $\Diamond$

\begin{coro}\label{coro23}
Let $p$, $w$, $n$, $M$, and $\vp$ be as in Theorem~\ref{theoQj_new1}. If $w^{-1}f\in L_{p}\cap L_{1}$ and
 $\omega_n(w^{-1}f,\delta)_{1}=\mathcal{O}(\delta^{d+a})$ as $\delta\to 0$, where $a\in (0,n-d)$, then
$$
\bigg\Vert f- m^{-j/2}\sum\limits_{k\in\zd} f(- M^{-j}k) \phi_{jk}\bigg\Vert_{p, 1/w}=\mathcal{O}\(|\lambda|^{-j(a+d/p)}\),\quad j\to \infty,
$$
where $\lambda$ is the smallest (in absolute value) diagonal element of $M$.
\end{coro}


{\bf Proof.} Applying the estimate
$$
|\mathcal{F}\(w^{-1}f\)(\xi)|\le C\omega_n\(w^{-1}f,{|\xi|^{-1}}\)_1=\mathcal{O}\(|\xi|^{-(d+a)}\),
$$
which can be found, e.g., in~\cite{Treb},
we see that all assumptions of Theorem~\ref{theoQj_new1} are satisfied.
It is also obvious that the required estimate holds for the second  term in the right hand side of~\eqref{8++++}.
To estimate the first term, we can use the well-known embedding for the Besov spaces (see, e.g.,~\cite[6.5.1 and 6.2.5]{BL}):
$$
B_{1,\infty}^{\alpha+d}(\R^d)\subset B_{p,\infty}^{a+d/p}(\R^d),
$$
which implies that $\omega_\nu(w^{-1}f,\delta)_p=\mathcal{O}(\delta^{a+d/p})$ for all integer $\nu\in (a+d/p,n]$. To estimate the moduli of smoothness of order $\nu\in [1, a+d/p]$, we can use the Marchaud inequality~\eqref{eqMarch} and the Jackson inequality given in Theorem~\ref{jackson}.
This provides the following estimate
$$
\sum_{\nu=0}^n \delta^{\nu q} \omega_{n-\nu}\(w^{-1}f,\delta\)_{p}^q=\mathcal{O}\(\delta^{q(a+d/p)}\),
$$
which together with Remark~\ref{rem1+++_} proves the corollary.
$\Diamond$

%

\end{document}